\documentclass[12pt,leqno]{article}
\usepackage{amsmath, amsthm, amsfonts, amssymb, color}
\usepackage{mathrsfs}
\usepackage{color}
\newtheorem{thm}{Theorem}[section]
\newtheorem{cor}[thm]{Corollary}
\newtheorem{lem}[thm]{Lemma}

\newtheorem{exa}[thm]{Example}
\theoremstyle{definition}
\newtheorem{defn}{Definition}[section]
\newcommand{\scr}[1]{\mathscr #1}
\definecolor{wco}{rgb}{0.5,0.2,0.3}

\numberwithin{equation}{section} \theoremstyle{remark}
\newtheorem{rem}{Remark}[section]

\title{{\bf The Onsager-Machlup action functional for degenerate McKean-Vlasov Stochastic Differential Equations}\footnote{The work is partially supported  by the NSFC No. 12171084 and the
		Fundamental Research Funds for the Central Universities No. RF1028623037.}
}

\author{
	{\bf   Shanqi Liu $^{a)}$ and Hongjun Gao $^{b)}$}\\
	\footnotesize{ a)School of Mathematical Science, Nanjing Normal University, Nanjing 210023, China}\\
	\footnotesize{ b)School of mathematics, Southeast University, Nanjing 211189, China}\\
	\footnotesize{  shanqiliumath@126.com,  hjgao@seu.edu.cn (Corresponding author)}}

\begin{document}

\def\R{\mathbb R}  \def\ff{\frac} \def\ss{\sqrt} \def\B{\mathbf
B}
\def\N{\mathbb N} \def\kk{\kappa} \def\m{{\bf m}}
\def\dd{\delta} \def\DD{\Delta} \def\vv{\varepsilon} \def\rr{\rho}
\def\<{\langle} \def\>{\rangle} \def\GG{\Gamma} \def\gg{\gamma}
  \def\nn{\nabla} \def\pp{\partial} \def\EE{\scr E}
\def\d{\text{\rm{d}}} \def\bb{\beta} \def\aa{\alpha} \def\D{\scr D}
  \def\si{\sigma} \def\ess{\text{\rm{ess}}}
\def\beg{\begin} 
\baselineskip 20pt
\maketitle
\begin{abstract} The purpose of this paper is to investigate the existence of the Onsager-Machlup action functional for degenerate McKean-Vlasov stochastic differential equations. To this end, we first derive Onsager-Machlup action functional for degenerate McKean-Vlasov stochastic differential equations with constant diffusion in a broad set of norms by Girsanov transformation, some conditioned exponential inequalities and It$\mathrm{\hat{o}}$ formulas for distribution dependent functional. Then an example is given to illustrate our results.

\end{abstract} \noindent

 2020 AMS subject Classification:\ 60H10; 37H10; 82C35.   \\
\noindent
 Keywords: Onsager-Machlup action functional;  degenerate stochastic differential equations; McKean-Vlasov; Girsanov transformation.\\
 \noindent
 \vskip 2cm

\section{Introduction}
The dynamical behavior of many systems arising in physics, chemistry, biology, etc, is often influenced by noise fluctuations. Among one of the interesting and challenging problems is to study transition theory between metastable states (or metastability) of stochastic systems. For the past several decades, this problem has encouraged many physical, mathematical and statistical studies in different backgrounds and motivations. We refer to \cite{GB, QD, DB, WE} and the references therein for an excellent review. In this paper, we focus on studying transition theory at the level of Onsager-Machlup (OM) action functional theory.

OM action functional was first given by Onsager and Machlup \cite{OM1, OM2} as the probability density functional for diffusion processes with linear drift and constant diffusion coefficients, and then Tizsa and Manning \cite{LT} generalized the results of Onsager and Machlup to nonlinear equations. The key point was to express the transition probability of a diffusion process by a functional integral over paths of the process, and the integrand was called OM action function. Then regarding OM action function as a Lagrangian,  the most probable path of the diffusion process was determined by variational principle. However, the paths of the diffusion process are almost surely nowhere differentiable. This means that it is not feasible to use the variational principle for the path of diffusion process. To modify that, Stratonovich \cite{RL} proposed the rigorous mathematical idea: one can ask for the probability that a path lies within a certain region, which may be a tube along a differentiable function (mostly called as reference path), comparing the probabilities of different tubes of the same 'thicknes', the OM action function is expressed by reference path instead of the path of diffusion process.

Indeed, classical Onsager-Machlup theory as shown in \cite{DB,NI}, for any reference path $\phi,\psi\in C^2([0,1], \mathbb{R}^m)$, the Onsager-Machlup action functional for $X_t$ is defined by
$$\lim_{\delta\to 0}\frac{P\big(\sup\limits_{t \in[0, T]}\big|X(t)-\phi(t)\big| \leq \delta\big)}{P\big(\sup\limits_{t \in[0, T]}\big|X(t)-\psi(t)\big| \leq \delta\big)}=\exp\Big(\frac{1}{\varepsilon^2}\big\{L_\varepsilon(\phi,\dot{\phi})-L_{\varepsilon}(\psi,\dot{\psi})\big\}\Big),$$
where $X(t)$ is the solution of non-degenerate stochastic differential equations (SDEs):
\begin{align}
	\d X(t)=f(t,X(t))\d t+\varepsilon\d W(t), X_0=x,\nonumber
\end{align}
where $f\in C_b^2([0,1]\times \mathbb{R}^{m})$, $W(t)$ is a $m$-dimensional Brownian motion and exact expression of Onsager-Machlup action functional is given by
$$L_\varepsilon(\phi,\dot{\phi})=-\frac{1}{2}\int_{0}^{1}\left|\dot{\phi}(t)-f(t,\phi(t))\right|^{2}\d t-\frac{\varepsilon}{2}\int_{0}^{1} \operatorname{div}_{x}f(t,\phi(t))\d t,$$	
where $\operatorname{div}_{x}$  denote the divergence on the $\phi(t)\in\mathbb{R}^m$. As noise intensity parameter $\varepsilon \to 0$,
Onsager-Machlup action functional coincides with Freidlin-Wentzell action functional \cite{FW} formally.

From the beginning of irreplaceable contributions of Stratonovich \cite{RL}, OM action functional theory was starting to receive considerable attention by  mathematicians. Many different approaches and new problems have arisen in the process. We first review works on OM action functional for stochastic differential equations driven by non-degenerate noise \cite{CG,MC1,MC2,YC,DB,AD,TF,GS,KH,KH1,NI,SM, RL,LS,YT}.  Ikeda and Watanabe \cite{NI} derived the OM action functional for reference path $\phi\in C^2([0,1],\mathbb{R}^{d})$ and taking the supremum norm $\|.\|_{\infty}$. D\"{u}rr and Bach \cite{DB} obtained the same results based on the Girsanov transformation of the quasi-translation invariant measure and the potential function (path integral representation). Shepp and Zeitouni \cite{LS} proved that this limit theorem still holds for every norm equivalent to the supremum norm and $\phi-x$ in the Cameron-Martin space. Capitaine \cite{MC1}  extended this result to a large class of natural norms on the Wiener space including particular cases of H\"{o}lder, Sobolev, and Besov norms. Hara and Takahashi provided in \cite{KH} a computation of the OM action functional of an elliptic diffusion process for the supremum norm and this result was extended in \cite{MC2} by Capitaine to norms that dominate supremum norm. In particular, the norms $\|\cdot\|$ could be any Euclidean norm dominating $L^2$-norm in the case of $\mathbb{R}^d$.

In addition to the derivation of OM action functional for SDEs derived by Brownian motion, the derivation of OM action functional for SDEs derived by fractional Brownian motion and L$\mathrm{\acute{e}}$vy process has also been explored to varying degrees. Nualart and Moret \cite{SM} first obtained the OM action functional in singular and regular cases respectively. In $\frac{1}{4}<H<\frac{1}{2}$ (singular case) $\|\cdot\|$ can be either the
supremum norm or $\alpha<H-\frac{1}{ 4}$ H\"{o}lder norm. In $H>\frac{1}{2}$ (regular case) the H\"{o}lder norm can only be taken as $H-\frac{1}{2}<\alpha< H-\frac{1}{4}$.  Recently in \cite{YC}, Chao and Duan first derived the OM action functional for scalar SDEs driven by the L$\mathrm{\acute{e}}$vy process with only small jumps by applying the idea of \cite{DB}. It is worth noting that the dimension considered in the above results is one-dimensional and the noise is additive.

At the same time, the intersection of OM action function theory with other theories and applications has also received a lot of attention, especially recently. Zeitouni and Dembo \cite{OZ, OZ1} first applied the idea of OM action functional to study the maximum a posterior (MAP) of trajectories for signal-observation system, which successfully establishes the relationship between OM action functional theory and nonlinear filtering theory. Then some developments in this direction see \cite{SA, BA1, MD}. More recently, $\Gamma$-convergence of Onsager-Machlup action functional \cite{QD, TL}, the bounds for the most probable transition time when time is not fixed \cite{YH2} and the connection between OM action functional and information projection problem \cite{ZS} are also considered. On the other hand, in specific application, climate change for global warming \cite{YZ}, single-species with Allee effect \cite{AT}, genetic regulatory network \cite{JH} and so on. OM action functional has become an important tool for understanding the dynamics of system.

However, we note that the noise of most of these results
is non-degenerate. But, so many mathematical models in finance, physics, biology etc (eg. \cite{QL, PW, ZW} and references therein), which the noise is not always non-degenerate. To characterize the most probable dynamics between the metastable states for degenerate SDEs in the sense of OM action functional theory, the most important problem is to determine the OM action functional for degenerate SDEs. This problem first attracted the attention of physicists. Kurchan \cite{JK} derived OM action functional via a Fokker-Planck equation \cite{HR} corresponding to the Langevin equation. Taniguchi and Cohen \cite{TT1, TT2} obtained the OM action functional for the Langevin equation by the path integral approach. Strict mathematical derivation originated from literature \cite{MCM} and \cite{SA} independently. Chaleyat-Maurel and Nualart \cite{MCM} derived Onsager-Machlup action functional for second-order stochastic differential equations with two-point boundary value condition. To derive the maximum likelihood state estimator for degenerate SDEs, Aihara and Bagchi \cite{SA} extended OM action functional into a degenerate version of OM action functional for reference path $\phi \in H^1([0,1],\mathbb{R}^d)$ with supremum norm by the approach of \cite{LS}.

McKean-Vlasov SDEs (also called distribution dependent SDEs or mean-field SDEs) describe limiting behaviors of individual particles in an interacting particle system of mean-field type when the number of particles goes to infinity (so called propagation of chaos) and widely used in various fields. One can refer to \cite{WL} and the references therein for more detials. Meanwhile, motivated by the study of OM action functional for non-degenerate McKean-Vlasov SDEs with additive noise \cite{SL}, it is natural to consider OM action functional for degenerate McKean-Vlasov SDEs with additive noise. To the best of our knowledge, this is the first result associated with OM action functional for degenerate McKean-Vlasov SDEs. The first difficulty is that how to connect different small probability of $\bar{X}^{1}_t$ and $\bar{X}^{2}_t$. The key point is using the special structure of the reference path. The second difficulty when we deal with degenerate McKean-Vlasov diffusion is that the distribution of $\bar{X}^{(2)}_t$ is independent of the limitation of the Brownian motion, which means that the distribution of $\bar{X}_t^{(2)}$ can not  be controlled in the general setting of small ball probability. So in order to overcame this difficulty, we would consider a class of degenerate McKean-Vlasov diffusion like $q(t,\bar{X}_t,\mathscr{L}_{\bar{X}^{(2)}_{t}})=\mu(\bar{X}_t)+\nu(\bar{X}_t)[E(\bar{X}^{2}_t)]^p\ (p\geq 1)$.

The remainder of this paper is organized as follows.  In Section $2$, we first recall some basic definitions, notations and lemmas for deriving the OM action function for degenerate McKean-Vlasov stochastic differential equations, and then we state our main results in Section $3$. In Section $4$, we apply our results to some models. Finally, complete proof of the main results are addressed
in Section $5$.
\section{Preliminaries and notation}\label{sec:opmr}

In this section, we recall some basic notations, assumptions and lemmas associated with degenerate McKean-Vlasov SDEs.
\subsection{Notations}
Throughout this paper, we denote $\mathbb{R}^{d}$ the $n$-dimensional Euclidean space, $\mathbb{R}^{d\times m}$ the collection of $d\times m$ matrices. And $I_{d\times d}$ and $O_{d\times d}$ denote identity matrix and null matrix respectively. Let $(\Omega,\mathcal{F},(\mathcal{F}_t)_{t\ge0},\mathbb{P})$ to be a given complete filtered probability space $(\Omega,\mathcal{F},(\mathcal{F}_t)_{t\ge0},\mathbb{P})$, where $\mathcal{F}_t$ is a nondecreasing family of sub-$\sigma$ fields of $\mathcal{F}$ satisfying the usual conditions. Let $\mathscr{P}$ be the space of all probability measures $\mu$ on $\mathbb{R}^{d}$, and let
$$\mathscr{P}_{2}(\mathbb{R}^{d})=\left\{\mu \in \mathscr{P}\left(\mathbb{R}^{d}\right): \mu\left(|\cdot|^{2}\right):=\int_{\mathbb{R}^{d}}|x|^{2} \mu(\mathrm{d} x)<\infty\right\}.$$
It is well known that $\mathscr{P}_{2}$ is a Polish space under the Wasserstein distance
$$
\mathbb{W}_{2}(\mu, \nu):=\inf _{\pi \in \mathscr{C}(\mu, \nu)}\left(\int_{\mathbb{R}^{m} \times \mathbb{R}^{m}}|x-y|^{2} \pi(\mathrm{d} x, \mathrm{d} y)\right)^{\frac{1}{2}}, \mu, \nu \in \mathscr{P}_{2}\left(\mathbb{R}^{d}\right),
$$
where $\mathscr{C}(\mu, \nu)$ is the set of couplings for $\mu$ and $\nu$; that is, $\pi \in \mathscr{C}(\mu, \nu)$ is a probability measure on $\mathbb{R}^{d} \times \mathbb{R}^{d}$ such that $\pi\left(\cdot \times \mathbb{R}^{d}\right)=\mu$ and $\pi\left(\mathbb{R}^{d} \times \cdot\right)=\nu$.
\begin{defn} \cite{PR}
	Let $t \in(0, 1]$.
	
	(i) A function $h: \mathscr{P}_{2}\left(\mathbb{R}^{d}\right) \rightarrow \mathbb{R}$ is called $L$-differentiable at $\mu \in \mathscr{P}_{2}\left(\mathbb{R}^{d}\right)$, if the functional
	$$
	L^{2}\left(\mathbb{R}^{d} \rightarrow \mathbb{R}^{d}, \mu\right) \ni \phi \mapsto h\left(\mu \circ(\operatorname{Id}+\phi)^{-1}\right)
	$$
	is Fréchet differentiable at $\phi=0 \in L^{2}\left(\mathbb{R}^{d} \rightarrow \mathbb{R}^{d}, \mu\right)$; that is, there exists (hence, unique) $\xi \in L^{2}\left(\mathbb{R}^{d} \rightarrow \mathbb{R}^{d}, \mu\right)$ such that
	$$
	\lim _{\mu\left(|\phi|^{2}\right) \rightarrow 0} \frac{h\left(\mu \circ(\mathrm{Id}+\phi)^{-1}\right)-h(\mu)-\mu(\langle\xi, \phi\rangle)}{\sqrt{\mu\left(|\phi|^{2}\right)}}=0 .
	$$
	In this case, we denote $\partial_{\mu} h(\mu)=\xi$ and call it the $L$-derivative of $h$ at $\mu$.
	
	(ii) A function $h: \mathscr{P}_{2}\left(\mathbb{R}^{d}\right) \rightarrow \mathbb{R}$ is called $L$-differentiable on $\mathscr{P}_{2}\left(\mathbb{R}^{d}\right)$ if the $L$-derivative $\partial_{\mu} h(\mu)$ exists for all $\mu \in \mathscr{P}_{2}\left(\mathbb{R}^{d}\right)$. If moreover $\left(\partial_{\mu} h(\mu)\right)(y)$ has a version differentiable in $y \in \mathbb{R}^{d}$ such that $\left(\partial_{\mu} h(\mu)\right)(y)$ and $\partial_{y}\left(\partial_{\mu} h(\mu)\right)(y)$ are jointly continuous in $(\mu, y) \in$ $\mathscr{P}_{2}\left(\mathbb{R}^{d}\right) \times \mathbb{R}^{d}$, we denote $h \in C^{(1,1)}\left(\mathscr{P}_{2}\left(\mathbb{R}^{d}\right)\right)$.
	
	(iii) A function $h:[0, 1] \times \mathbb{R}^{d} \times \mathscr{P}_{2}\left(\mathbb{R}^{d}\right) \rightarrow \mathbb{R}$ is said to be in the class $C^{1,2,(1,1)}([0, 1] \times$ $\mathbb{R}^{d} \times \mathscr{P}_{2}\left(\mathbb{R}^{d}\right))$, if the derivatives
	$$
	\partial_{t} h(t, x, \mu), \partial_{x} h(t, x, \mu), \partial_{x}^{2} h(t, x, \mu), \partial_{\mu} h(t, x, \mu)(y), \partial_{y} \partial_{\mu} h(t, x, \mu)(y)
	$$
	exist and are jointly continuous in the corresponding arguments $(t, x, \mu)$ or $(t, x, \mu, y)$. If $f \in C^{1,2,(1,1)}([0, 1] \times \mathbb{R}^{d}\times \mathscr{P}_{2}\left(\mathbb{R}^{d}\right))$ with all these derivatives bounded on $[0, 1] \times$ $\mathbb{R}^{d} \times \mathscr{P}_{2}\left(\mathbb{R}^{d}\right)$, we denote $f \in C_{b}^{1,2,(1,1)}([0, 1] \times \mathbb{R}^{d} \times \mathscr{P}_{2}\left(\mathbb{R}^{d}\right))$.
	
	(iv) Finally, we write $h \in \mathscr{C}\left([0, 1] \times \mathbb{R}^{d} \times \mathscr{P}_{2}\left(\mathbb{R}^{d}\right)\right)$, if $h \in C^{1,2,(1,1)}\left([0, 1] \times \mathbb{R}^{d} \times \mathscr{P}_{2}\left(\mathbb{R}^{d}\right)\right)$ and the function
	$$
	(t, x, \mu) \mapsto \int_{\mathbb{R}^{d}}\left\{\left\|\partial_{y} \partial_{\mu} h\right\|+\left\|\partial_{\mu} h\right\|^{2}\right\}(t, x, \mu)(y) \mu(\mathrm{d} y)
	$$
	is locally bounded, i.e. it is bounded on compact subsets of $[0, 1] \times \mathbb{R}^{d} \times \mathscr{P}_{2}\left(\mathbb{R}^{d}\right)$.
\end{defn}
\subsection{Degenerate McKean-Vlasov SDEs}	
Consider the following McKean-Vlasov SDE  on $\mathbb{R}^{d+m}$:
\begin{align}\label{all}
	\d Z_t=b(t,Z_t,\mathscr{L}_{Z_{t}})\d t+\sigma\d \tilde{W}(t), Z_0=z,
\end{align}
where $b: [0,1]\times \mathbb{R}^{d+m} \times \mathscr{P}_{2}\left(\mathbb{R}^{d+m}\right) \rightarrow \mathbb{R}^{d+m}$ is measurable, $\sigma$ is a $d+m\times d+m$ matrix, $\tilde{W}(t)$ is a $d+m$-dimensional Brownian motion. $\mathscr{L}_{Z_{t}}$ is the law of $Z_{t}$ under the given completed filtration probability
space $(\Omega,\mathcal{F},(\mathcal{F}_t)_{t\ge0},\mathbb{P})$. In particular, if $\phi(t)$ is a                                                                                                                                                                                                                                                                                                                                                                                                                                                                                                                                                                                                                                                                                                                                                                                                                                                                                                                                                                                                                                                                                                                                                                                                                                                                                                                                                                                                                                                                                                                                                                                                                                                                                                                                                                                                                                                                                                                                                                                                                                                                                                                                                                                                                                                                                                                                                                                                                                                                                                                                                                                                                                                                                                                                                                                                                                                                                                                                                                                                                                                                                                                                                                                                                                                                                                                                                                                                                                                                                                                                                                                                                                                                                                                                                                                                                                                                                                                                                                                                                                                                                                                                                                                                                                                                                                                                                                                                                                                                                                                                                                                                                                                                                                                                                                                                                                                                                                                                                                                                                                                                                                                                                                                                                                                                                                                                                                                                                                                                                                                                                                                                                                                                                                                                                                                                                                                                                                                                                                                                                                                                                                                                                                                                                                                                                                                                                                                                                                                                                                                                                                                                                                                                                                                                                                                                                                                                                                                                                                                                                                                                                                                                                                                                                                                                                                                                                                                                                                                                                                                                                                                                                                                                                                                                                                                                                                                                                                                                                                                                                                                                                                                                                                                                                                                                                                                                                                                                                                                                                                                                                                                                                                                                                                                                                                                                                                                                                                                                                                                                                                                                                                                                                                                                                                                                                                                                                                                                                                                                                                                                                                                                                                                                                                                                                                                                                                                                                                                                                                                                                                                                                                                                                                                                                                                                                                                                                                                                                                                                                                                                                                                                                                                                                                                                                                                                                                                                                                                                                                                                                                                                                                                                                                                                                                                                                                                                                                                                                                                                                                                                                                                                                                                                                                                                                                                                                                                                                                                                                                                                                                                                                                                                                                                                                                                                                                                                                                                                                                                                                                                                                                                                                                                                                                                                                                                                                                                                                                                                                                               deterministic path, then the law of $\phi(t)$ is a Dirac measure at $\phi(t)$, i.e. $\mathscr{L}_{\phi(t)}=\delta_{\phi(t)}$.

For the case of $\sigma=I_{{d+m}\times {d+m}}$.  In \cite{SL}, we proved that if we assume $b\in C_b^{1,2,(1,1)}([0, 1] \times$ $\mathbb{R}^{d} \times \mathscr{P}_{2}\left(\mathbb{R}^{d}\right))$ and  reference path $\phi-x$ belongs to Cameron-Martin $\mathcal{H}$. Then under some suitable conditions, the OM action functional of $Z_t$ for any norm dominating $L^2([0,1],\mathbb{R}^{d+m})$  exists and is given by:
$$L(t,\phi,\dot{\phi},\delta_{\phi})=-\frac{1}{2}\int_{0}^{1}|\dot{\phi}(t)-b(t,\phi(t),\mathscr{L}_{\phi(t)})|^2\ \d t-\frac{1}{2}\int_{0}^{1}\operatorname{div}_{x}b(\phi(t),\mathscr{L}_{\phi(t)})\d t,$$
where $\operatorname{div}_{x}$  denote the divergence on the $\phi(t)\in\mathbb{R}^{d+m}$. And one can easily remark that when the drift coefficient $b$ does not depend on the distribution of the solution, i.e. $b(t,Z_t,\mathscr{L}_{Z_{t}})=b(t,Z_t)$, OM action functional of \eqref{all} coincides with the classical OM action functional.

In this paper, we shall consider \eqref{all} with degenerate noise.
Let $$\sigma= \Xi:=\left(\begin{array}{cc}
	O_{d\times d} & O_{d\times m}\\
	O_{m\times d} & I_{m\times m}
\end{array}\right),$$
and
$$Z_t:=\bar{X}_t=\left(\begin{array}{c}
	\bar{X}^{(1)}_t \\
	\bar{X}^{(2)}_t
\end{array}\right), \bar{X}_0=x^{\mu}:=\left(\begin{array}{c}
	x^{(1)} \\
	x^{(2)}
\end{array}\right), b(\bar{X}_t, \mathscr{L}_{\bar{X}_{t}}):=\left(\begin{array}{c}
	p\left(\bar{X}_{t}\right) \\
	q\left(\bar{X}_{t}, \mathscr{L}_{\bar{X}^{(2)}_{t}}\right)
\end{array}\right).$$
Then we have the following McKean-Vlasov stochastic Hamiltonian system on $\mathbb{R}^{d} \times \mathbb{R}^{m}$ :
\begin{equation}\label{eq1}
	\left\{\begin{array}{l}
		\mathrm{d} \bar{X}_{t}^{(1)}=p\left(\bar{X}_{t}\right) \mathrm{d} t,\\
		\mathrm{d} \bar{X}_{t}^{(2)}=q\left(\bar{X}_{t}, \mathscr{L}_{\bar{X}^{(2)}_{t}}\right) \mathrm{d} t+ \mathrm{d} W(t),
	\end{array}\right.
\end{equation}
where the coefficients of diffusion term and drift term in \eqref{eq1}
$$
p: \mathbb{R}^{d+m}\rightarrow \mathbb{R}^{d}, q: \mathbb{R}^{d+m} \times \mathscr{P}_{2}\left(\mathbb{R}^{m}\right) \rightarrow \mathbb{R}^{m}
$$
are measurable and $\mathscr{P}_{2}(\cdot)$ is defined as above. The well-posedness for \eqref{eq1}, regularity assumptions about $p, q$ will be given later.

\subsection{Onsager-Machlup action functional}	
The problem we are interested in is to derive the Onsager-Machlup action functional for stochastic process $Z_t$ between two metastable states, that is, how to determine the maximal probability of the original path $Z_t$ surrounding the reference path $\phi(t)$:
$$	\mathbb{P}(\|Z-\phi\|<\varepsilon),$$
where $\|\cdot\|$ denotes the norm dominates $L^2[0,1]$, and reference path $\phi(t)$ can be regarded as the smooth enough path to approximate the original path. By estimating the probability of the original path surrounding the reference path, the Onsager-Machlup Lagrangian is the quantity that does not change with the thickness of the tube.  Finally, it is worth noting that we  use $\|\cdot\|$ to denote the norm that dominates  $L^2([0,1],\mathbb{R}^{k})$ with $k=d,m$ or $d+m$.

We now state the definitions of Onsager-Machlup action functional and Onsager-Machlup action function.
\begin{defn}
	Consider a tube surrounding a reference path $\phi(t)$, if for $\varepsilon$ sufficiently small we estimate the probability of the solution process $Z_t$ lying in this tube in the form:
	$$\mathbb{P}(\|Z-\phi\|\le\varepsilon)\propto C(\varepsilon)\exp\Big\{\int_{0}^{1}\mathnormal{OM}(\phi,\dot{\phi},\delta_{\phi})\Big\},$$
	then integrand $\mathnormal{OM}(\phi,\dot{\phi},\delta_{\phi})$ is called  Onsager-Machlup action function. Where $\propto$ denotes the equivalence relation for $\varepsilon$ small enough and $\|\cdot\|$ is a suitable norm. We also call $\int_{0}^{1}\mathnormal{OM}(\phi,\dot{\phi},\delta_{\phi})\mathrm{d}t$ the  Onsager-Machlup action functional.
\end{defn}
\subsection{The structure of reference path}
The effect of introducing the reference path is not only to overcome the regularity that the original path does not have, for stochastic differential equations with a specific structure, sometimes we also need to formulate a similar structure for the reference path. For example, recall previous works, Bardina, Rovira and Tindel \cite{XB1} constructed the structure of reference path based on the structure of the original stochastic evolution equation and the choice of mild solution; Moret and Nualart in \cite{SM} constructed reference path based on the connection between Brownian motion and fractional Brownian motion. And Chaleyat-Maurel and Nualart in \cite{MCM} constructed reference path based on the structure of second-order SDEs and two-point boundary value conditions. Hu \cite{YH1} constructed reference path based on the connection between Brownian motion and geometric Brownian motion with logarithmic transformation. All these imply that it is reasonable to construct an appropriate structure for the reference path when dealing with complex situations. This idea is also applicable to derive OM action functional for degenerate stochastic differential equations. The structure of the reference path for degenerate SDEs first given by \cite{SA}. In our paper we also need the reference path to keep the coupling effect:
$$\phi(t)-x^{\mu}=(\phi^1(t)-x^{(1)},\phi^2(t)-x^{(2)})\in \text{Cameron-Martin space } \mathcal{H}\quad \text{and} \quad \d\phi^1(t)=p(\phi(t), \mathscr{L}_{\phi(t)})\d t,$$
where $\phi^1(t)\in\mathbb{R}^{d}, \phi^2(t)\in\mathbb{R}^{m}$ and $\phi^1(0)=x^{(1)}\in\mathbb{R}^d,\phi^2(t)=x^{(2)}\in\mathbb{R}^m$. Cameron-Martin space $\mathcal{H}$  stands for  the class  of all absolutely continuous functions $h$
such that $h(0) = 0$ and $\dot{h}\in L^2[0,1]$, the inner product is given by
the formula
$$(h_1, h_2)_{\mathcal{H}}:=\int_{0}^{1}\dot{h}_1(t)\dot{h}_2(t)\d t.$$
\subsection{Assumptions}
Next, we impose some assumptions throughout this paper. And the norm $\|\cdot\|$ stand for a class of norms that dominate $L^2[0,1]$ and satisfy the following assumptions $\mathbf{(H_0-H_4)}$.

$\bullet$ $\mathbf{(H_0)}$
Assume there exists an increasing function $K:[0,1]\mapsto (0,\infty)$ such that

(i) $$|b(t,x,\mu)-b(t,y,\nu)|\leq K(t)\Big(|x-y|+\mathbb{W}_{2}(\mu,\nu)\Big),$$
for $t\in[0,1]$, $x,y\in\mathbb{R}^{d+m},\nu,\mu\in \mathscr{P}_{2}\left(\mathbb{R}^{d+m}\right)$.

(ii)
$$|b(t,0,\delta_0))|\leq K(t),$$
where $t\in [0,1], 0\in \mathbb{R}^{d+m}$ and $\delta_0$ is the Dirac measure at 0.

(iii) Assume reference path $\phi(t)-x^{\mu}=(\phi^1(t)-x^{1},\phi^2(t)-x^{2})\in \text{Cameron-Martin space } \mathcal{H}\quad \text{and} \quad \d\phi^1(t)$\\$=p(\phi(t))\d t$, then for every $c\in\mathbb{R}$,
\small{\begin{align}
		\limsup _{\varepsilon \rightarrow 0} E\Big[\exp\Big[c\Big(\int_{0}^{1}\Big|q(\tilde{X}^{(1)}_t,\phi^2(t)+W(t),\mathscr{L}_{\phi^{2}(t)+W(t)})-q(\phi^1(t),\phi^2(t),\mathscr{L}_{\phi^2(t)})\Big|^2\d t\Big)^{\frac{1}{2}}\Big]\vert \|W\|\leq \varepsilon\Big]\leqslant 1,\nonumber
\end{align}}
\begin{align}
	\limsup _{\varepsilon \rightarrow 0} E\Big[\exp\Big[ c\Big(\int_{0}^{1}\Big|\partial_{x^2}q(\tilde{X}^{(1)}_t, \phi^2(t)+W(t),\mathscr{L}_{\phi^2(t)+W(t)})-\partial_{x^2}q(\phi^1(t),\phi^2(t), \mathscr{L}_{\phi^{2}(t)})\Big|^2\d t\Big)^{\frac{1}{2}}\Big]\vert \|W\|\leq \varepsilon\Big]\leqslant 1,\nonumber
\end{align}
where $\tilde{X}^{(1)}_t$ is defined in \eqref{X^1} and
\begin{align}
	&\int_{0}^{1}\Big|W_i(t)\Big(\partial_{x_j^2}q(\tilde{X}^{(1)}(t), \phi^2(t)+W(t),\mathscr{L}_{\phi^2(t)+W(t)})-\partial_{x_j^2}q(\phi^1(t),\phi^2(t), \mathscr{L}_{\phi^{2}(t)})\Big)\Big|^2\d
	t\nonumber\\&\leq K_{8}\Big[\int_{0}^{1}|\tilde{X}^{(1)}_t|^4\d t+\int_{0}^{1}|W(t)|^4\d t\Big],\nonumber
\end{align}
for $\{\|W\|\leq\varepsilon\}$, $1\leq i,j \leq d$ and $K_{8}$ is positive constant.

$\bullet$ $\mathbf{(H_1)}$ $\text{(Geometrical properties on $\|\cdot\|$)}$ The norm $\|\cdot\|$ is invariant under the action of the orthogonal group $\mathcal{O}(\mathbb{R}^m)$ on the coordinates of the Brownian motion. i.e. for every $1\leq i\leq m$,
$$\|(W_1(t),\cdots, W_i(t),\cdots, W_m(t))\|=\|(W_1(t),\cdots, -W_i(t),\cdots, W_m(t))\|.$$

$\bullet$ $\mathbf{(H_2)}$ For every $1\leq i\leq m$ and every $c\in \mathbb{R}$,
\begin{align}
	\limsup _{\varepsilon \rightarrow 0} E[\exp(c|W_i(1)|^2)\vert \|W\|\leq \varepsilon]\leqslant 1.\nonumber
\end{align}

$\bullet$ $\mathbf{(H_3)}$
$\text{(Small ball probability estimation)}$ There exists $0<q<p\wedge 4$ such that for any $\varepsilon$ small enough
$$\int_{0}^{1}|W(t)|^4\d t\leq C_2\varepsilon^{p},$$
for $\|W\|\leq\varepsilon$ and
$$P(\|W\|\leq\varepsilon)\geq\exp(-\frac{C_3}{\varepsilon^q}),$$
where $C_2, C_3$ are positive constants.

$\bullet$  $\mathbf{(H_4)}$ $\text{(Controllability on small ball probability)}$ There exists a constant $C_5$ such that a.s. $\omega\in \Omega$,
$$\|\bar{X}^{(1)}-\phi^1\|\leq C_5 \|\bar{X}^{(2)}-\phi^2\|.$$

Here, we will systematically explain the motivation and purpose of these assumptions in the following remarks.
\begin{rem}
	As is well known, by $\mathbf{(H_0)}$ $(i), (ii)$, for any $s\geq 0$ and determinitic initial value $z$, \eqref{all} has a unique solution $(Z_t)_{t\geq s}$ with
	$$E\Big(\sup_{t\in [s, 1]}|Z_t|^2\Big)<\infty,$$
	and the imposed condition $(iii)$ in $\mathbf{(H_0)}$ is to extract divergence part of OM action functional for degenerate McKean-Vlasov SDEs \eqref{eq1}. Next we would verify $\mathbf{(H1)}$ $(iii)$ by some examples. Assume $q$ can be divided into $\mu(\tilde{X}_t)+v(\tilde{X}_t)[E(\tilde{X}^{(2)}_t)]^p\ (p\geq 1)$ and $q\in C_{b}^{1,2,(1,1)}\left([0, 1] \times \mathbb{R}^{\d+m} \times \mathscr{P}_{2}\left(\mathbb{R}^{m}\right)\right)$ in the following, which implies that $\mu,v\in C^2_{b}([0,1]\times \mathbb{R}^{\d+m})$.  So by using the fact $\mu,v$ are Lipschitz continuous and bounded with Lipschitz constant $K_1, K_2$ respectively, we have that
	\begin{align}
		&c\Big(\int_{0}^{1} \Big|q(\tilde{X}^{(1)}_t, \phi^2(t)+W(t),\mathscr{L}_{\phi^2(t)+W(t)})-q(\phi^{(1)}(t),\phi^{(2)}(t),\mathscr{L}_{\phi^{(2)}(t)})\Big|^2 \mathrm{d}t\Big)^{\frac{1}{2}}\nonumber\\&=c\Big(\int_{0}^{1}\Big| \mu(\tilde{X}^{(1)}_t, \phi^2(t)+W(t))+v(\tilde{X}^{(1)}_t, \phi^2(t)+W(t))[E( \phi^2(t)+W(t))]^p\nonumber\\&-\mu(\phi^1(t),\phi^2(t))-v(\phi^1(t),\phi^2(t))[E(\phi^2(t))]^p\Big|^2 \mathrm{d}t\Big)^{\frac{1}{2}}\nonumber\\&\leq \Big[\int_{0}^{1}3c^2\Big(K^2_1\big[|\tilde{X}_t^{(1)}-\phi^1(t)|^2+|W(t)|^2\big]+K^2_2\big[|\tilde{X}_t^{(1)}-\phi^1(t)|^2+|W(t)|^2\big](E(\phi^2(t)+W(t)))^{2p}\nonumber\\&+v^2(\phi^1(t),\phi^2(t))\big([E(\phi^2(t)+W(t))]^p-[E(\phi^2(t))]^p\big)^2\Big)\d t\Big]^{\frac{1}{2}}\nonumber\\&\leq
		K_3\Big(\int_{0}^{1}|W(t)|^2\d t+\|W\|^2\Big)^{\frac{1}{2}}\leq 2K_3\|W\|\leq 2K_3 \varepsilon,\nonumber
	\end{align}
	where the last step Minkowski's inequality and assumption $\mathbf{(H^*_4)}$ are used. Similarly we have
	\begin{align}
		&c\Big(\int_{0}^{1}\Big|\partial_{x^2}q(\tilde{X}^{(1)}_t, \phi^2(t)+W(t),\mathscr{L}_{\phi^2(t)+W(t)})-\partial_{x^2}q(\phi^1(t),\phi^2(t), \mathscr{L}_{\phi^{2}(t)})\Big|^2\d t\Big)^{\frac{1}{2}}\nonumber\\&\leq
		K_4\Big(\int_{0}^{1}|W(t)|^2\d t+\|W\|^2\Big)^{\frac{1}{2}}\leq 2K_4\|B\|\leq 2K_4 \varepsilon,\nonumber
	\end{align}
	and
	\begin{align}
		&\int_{0}^{1}\Big|W_i(t)\Big(\partial_{x_j^2}q(\tilde{X}^{(1)}(t), \phi^2(t)+W(t),\mathscr{L}_{\phi^2(t)+W(t)})-\partial_{x_j^2}q(\phi^1(t),\phi^2(t), \mathscr{L}_{\phi^{2}(t)})\Big)\Big|^2\d
		t\nonumber\\&\leq K_{8}\Big[\int_{0}^{1}|\tilde{X}^{(1)}_t|^4\d t+\int_{0}^{1}|W(t)|^4\d t\Big],\nonumber
	\end{align}
	for $\{\|W\|\leq\varepsilon\}$, $1\leq i,j \leq d$ and $K_{8}$ is positive constant.
	Therefore $\mathbf{(H_0)}$ $(iii)$ is verfied.
\end{rem}
\begin{rem}
	The assumption $\mathbf{(H_1)}$ first appears in \cite{LS1}, and they called the norm is convex norm if it satisfies this condition. In \cite{MC1}, Capitaine gave some equivalent forms. The purpose of introducing this condition is to satisfy the condition of Lemma \ref{Cross lemma}.
\end{rem}
\begin{rem}
	For the assumption $\mathbf{(H_2)}$, it can be found in (P2) of \cite{MC1}, the purpose is to deal with the endpoint value after applying It$\mathrm{\hat{o}}$ formula or Lemma \ref{aa} for stochastic integral $\int_{0}^{1}a(t)W_i(t)\d W_i(t)$, where $a(t)$ is a bounded function.
\end{rem}
\begin{rem}
	Assumption $\mathbf{(H_3)}$ is a usual assumption to control remainder term by small ball probability, we can note similar assumptions in previous works. For example, in the case of finite-dimensional Brownian motion, Capitaine \cite{MC1} gave assumptions $(P1)$ and $(P2)$. In the case of infinite-dimensional cylindrical Brownian motion, Bardina, Rovira and Tindel \cite{XB1} imposed assumption $(\bar{h}1)$ and $(\bar{h}2)$. In the case of finite-dimensional fractional Brownian motion, Theorem $7$ and Theorem $8$ in \cite{SM} implied small ball probability assumption (the limit of the Hurst index $H$ is mainly influenced by the limit of the small ball probability estimation). As far as we know, small ball probability estimation is now the most classical tool for dealing with remainder terms, this idea first originated in \cite{NI}. Inspired by these works, we modified them to satisfy our needs when dealing with degenerate situations. More specifically, we replace $0<q<p$ by  $0<q<p\wedge 4$.
\end{rem}
\begin{rem}\label{equiva of small ball}
	Finally, we would explain the controllability on small ball probability. When deriving the OM action functional in the degenerate case, an important idea is that we should establish a connection between the two component solution processes, as well as reference path. Therefore, in the first work \cite{SA} of  deriving OM action functional for degenerate stochastic differential equations,  this controllability on small ball probability is obvious for supremum norm by Gronwall inequality. But since in our paper the norm $\|\cdot\|$ stands for a class of abstract norms that dominate $L^2[0,1]$, we can't prove $\mathbf{(H_4)}$ directly. So in order to verify the rationality of assumptions, we have verified the different norms  such as uniform norm, $L^p$-norms and H$\ddot{o}$lder-norm  in the Appendix. Furthermore,
	$\mathbf{(H_4)}$ is equivalent to the following form:
	$$\mathbb{P}(\|\bar{X}-\phi\|\leq \varepsilon)=\mathbb{P}(\|\bar{X}^{(2)}-\phi^2\|\leq \varepsilon),$$
	for any norms $\|\cdot\|$ dominating $L^2[0,1]$.
\end{rem}
\subsection{Technical lemmas}

Before we derive Onsager-Machlup action functional for degenerate SDEs, some technical lemmas are needed.
\begin{lem}[\cite{NI} \text{pp 536-537}]\label{Separation lemma}
	For a fixed $n\geqslant1$, let $I_1, \ldots, I_n $ be $n$ random variables defined on $(\Omega,\mathcal{F},\mathbb{P})$ and $\{A_\varepsilon;\varepsilon >0\}$ a family of sets in $\mathcal{F}$. Suppose that for any $c\in \mathbb{R}$ and any $i=1,\dots, n$, if we have
	\begin{align}
		\limsup _{\varepsilon \rightarrow 0} E[\exp(cI_i)\vert A_\varepsilon]\leqslant 1.\nonumber
	\end{align}
	Then
	\begin{align}
		\lim_{\varepsilon\to 0} E \left [\exp \left(\sum_{i=1}^{n}cI_i \right ) \bigg| A_\varepsilon \right ]= 1.\nonumber
	\end{align}
\end{lem}
\begin{lem}\cite{LS}\label{no random lemma}
	Let $f$ be a deterministic function in $L^2[0,1]$. Define $I_i(f)=\int_{0}^{1}f(t)\d W_i(t)$. If the norm $\|\cdot\|$ dominates the $L^1[0,1]$-norm then
	\begin{align}
		\lim_{\varepsilon\to 0} E [\exp(|I_i(f)|)|\|				W\|<\varepsilon]= 1.\nonumber
	\end{align}
\end{lem}
\begin{lem}[\cite{LS1} \text{Theorem 1}]\label{Cross lemma}
	Assume that the norm $\|\cdot\|$ dominates the $L^{1}[0,1]$-norm and satisfies
	for every $1\leq i\leq m$,
	$$\|(W_1(t),\cdots, W_i(t),\cdots, W_m(t))\|=\|(W_1(t),\cdots, -W_i(t),\cdots, W_m(t))\|.$$
	Let $\mathscr{F}_{i}$ be the $\sigma$-field generated by $\left\{W_1(t), \ldots, W_{i-1}(t), W_{i+1}(t), \ldots, W_m(t) ; 0 \leq\right.$ $t \leq 1\} .$ Let $\Psi(\cdot)$ be an $\mathscr{F}_{i}$-adapted function such that
	$$
	\forall c \in \mathbb{R}^{+}, \lim _{\epsilon \rightarrow 0} E\left(\exp \left\{c \int_{0}^{1} \Psi^{2}(t) \d t\right\} \mid\|W\|<\varepsilon\right)=1.
	$$
	Then
	$$
	\lim _{\epsilon \rightarrow 0} E\left(\exp \left\{\left|\int_{0}^{1} \Psi(t) \d W_i(t)\right|\right\} \mid\|W\|<\varepsilon\right)=1.
	$$
\end{lem}

\begin{lem}[\cite{PR} \text{Lemma 3.1}](It$\hat{o}$ formulas  for distribution dependent functional)\label{aa}
	Assume $Z_t$ satisfies $(2.1)$, then for any $h \in \mathscr{C}([0, 1] \times$ $\left.\mathbb{R}^{d+m} \times \mathscr{P}_{2}\left(\mathbb{R}^{d+m}\right)\right), h\left(t, Z_{t}, \mathscr{L}_{Z_{t}}\right)$ is a semi-martingale with
	$$
	\mathrm{d} h\left(t, Z_{t}, \mathscr{L}_{Z_{t}}\right)=\left(\partial_{t}+\mathbf{L}_{b, \sigma}\right) h\left(t, Z_{t}, \mathscr{L}_{Z_{t}}\right) \mathrm{d} t+\left\langle\left(\sigma^{*} \partial_{x} h\right)\left(t, Z_{t}, \mathscr{L}_{Z_{t}}\right), \mathrm{d}\tilde{W}(t)\right\rangle,
	$$
	where
	$$
	\begin{aligned}
		\mathbf{L}_{b, \sigma} h(t, x, \mu)=& \frac{1}{2} \operatorname{tr}\left(\sigma \sigma^{*} \partial_{x}^{2} h\right)(t, x, \mu)+\left\langle b, \partial_{x} h\right\rangle(t, x, \mu) \\
		&+\int_{\mathbb{R}^{d+m}}\Big[\frac{1}{2} \operatorname{tr}\left\{\left(\sigma \sigma^{*}\right)(t, y, \mu) \partial_{y} \partial_{\mu} h(t, x, \mu)(y)\right\}+\langle b(t,y,\mu), \partial_{\mu}h(t,x,\mu)(y)\rangle\Big]\mu(\d y).
	\end{aligned}
	$$
\end{lem}
Next, we introduce the Moore-Penrose generalized inverse, which will be used to characterize the result of Theorem $3.1$ and Theorem $3.6$.
\begin{defn}\cite{AB}
	The Moore-Penrose generalized inverse matrix of  $A\in\mathbb{R}^{d\times n}$, denoted by $A^{+}$, is an $n\times d$ matrix if  the following four defining properties are satisfied:
	$$
	\begin{aligned}
		\left(A A^{+}\right)^{*} &=A A^{+}, & &\left(A^{+} A\right)^{*}=A^{+} A, \\
		A A^{+} A &=A, & & A^{+} A A^{+}=A^{+}.
	\end{aligned}
	$$
	It is unique if it exists.
\end{defn}
We also have explicit formula for the Moore-Penrose inverse of an $m\times n$ partitioned matrix $M=\left(\begin{array}{cc}A & D \\ B & C\end{array}\right)$.
\begin{lem}\cite{HH}\label{MP inverse}
	Let $M=\left(\begin{array}{cc}A & D \\ B & C\end{array}\right)$. Then
	$$
	M^{+}=\left(\begin{array}{cc}
		K^{+}\left(A^{*}-E F\right) & K^{+}\left(B^{*}-E H\right) \\
		F & H
	\end{array}\right),
	$$
	where
	$$
	\begin{aligned}
		&K=A^{*} A+B^{*} B, \\
		&E=A^{*} D+B^{*} C, \\
		&R=D-A K^{+} E, \\
		&S=C-B K^{+} E, \\
		&L=R^{*} R+S^{*} S, \\
		&T=K^{+} E\left(I-L^{+} L\right), \\
		&F=L^{+} R^{*}+\left(I-L^{+} L\right)\left(I+T^{*} T\right)^{-1}\left(K^{+} E\right)^{*} K^{+}\left(A^{*}-E L^{+} R^{*}\right),
	\end{aligned}
	$$
	and
	$$
	H=L^{+} S^{*}+\left(I-L^{+} L\right)\left(I+T^{*} T\right)^{-1}\left(K^{+} E\right)^{*} K^{+}\left(B^{*}-E L^{+} S^{*}\right).
	$$
\end{lem}
We give two simple examples to help us understand below.
\begin{exa}
	If we consider the Moore-Penrose generalized inverse matrix of $A=(1,1)$, by the defining properties we can easily get
	$$A^{+}=\left(\begin{array}{cc}\frac{1}{2} \\ \frac{1}{2} \end{array}\right).$$
\end{exa}
\begin{exa}
	Here we consider a very simple partitioned matrix have introduced in Case $I$, the Moore-Penrose generalized inverse matrix of $$\Xi=\left(\begin{array}{cc}
		O_{d\times d} & O_{d\times m} \\
		O_{m\times d} & I_{m\times m}
	\end{array}\right),$$
	by Lemma $2.13$, we can get
	$$\Xi^{+}=\left(\begin{array}{cc}
		O_{d\times d} & O_{d\times m} \\
		O_{m\times d} & I_{m\times m}
	\end{array}\right),$$
	clearly, $\Xi=\Xi^{+}$ and $\Xi\Xi^{+}\neq I_{(m+d)\times (m+d)}$.
\end{exa}
\section{Main Results} 
The purpose of our paper is to study OM action functional theory for degenerate McKean-Vlasov SDEs. We now state results associated with OM action functional for degenerate McKean-Vlasov SDEs.
\begin{thm} \textbf{(The Onsager-Machlup action functional for degenerate McKean-Vlasov SDEs with additive noise)}\label{thm 2}
	For a class of degenerate McKean-Vlasov stochastic differential equation in the form of \eqref{eq1}, assume that $p \in C_{b}^{1}\left(\mathbb{R}^{d+m}\right)$, $q \in C_{b}^{2,(1,1)}\left(\mathbb{R}^{d+m} \times \mathscr{P}_{2}\left(\mathbb{R}^{m}\right)\right)$ and  assumptions $\mathbf{(H_0)}, \mathbf{(H_1)}, \mathbf{(H_2)},\mathbf{(H_3)}, \mathbf{(H_4)}$ hold. Reference path $\phi(t)-x^\mu$ belongs to Cameron-Martin $\mathcal{H}$. Then the Onsager-Machlup action functional of $\bar{X}_t$ for any norm dominating $L^2[0,1]$  exists and is given by
	$$L(\phi,\dot{\phi},\delta_{\phi^2})=-\frac{1}{2}\int_{0}^{1}\Big|\dot{\phi}^2(t)-q(\phi(t),\delta_{\phi^2(t)})\Big|^2\d t-\frac{1}{2}\int_{0}^{1}\operatorname{div}_{x^2}q(\phi(t),\delta_{\phi^2(t)})\d t,$$
	where $\operatorname{div}_{x^2}$  denote the divergence on the $\phi^2(t)\in\mathbb{R}^m$, or the following equivalent global expression:
	$$
	L(\phi, \dot{\phi},\delta_{\phi^2})=-\frac{1}{2}\int_{0}^{1}\Big|\Xi^{+}[\dot{\phi}(t)-b(\phi(t), \delta_{\phi^2(t)})]\Big|^2\d t-\frac{1}{2}\int_{0}^{1} \operatorname{div_x}[\Xi^{+}\cdot b](\phi(t), \delta_{\phi^2(t)})\d t,
	$$
	where $\operatorname{div_x}$ denote the divergence on the $\mathbb{R}^{d+m}$ and $\phi$ is the function such that
	$$\phi(t)-x^{\mu}=(\phi^1(t)-x^{(1)},\phi^2(t)-x^{(2)})\in \text{Cameron-Martin space } \mathcal{H}\quad \text{and} \quad \d\phi^1(t)=p(\phi(t))\d t.$$
\end{thm}
\begin{cor}
	Under the same assumptions of Theorem \ref{thm 2}, in particular, if we let $p\left(\bar{X}_{t}\right) =g(\bar{X}^{(1)}_{t})$, then it is clear $\bar{X}^{(1)}_t$ is a deterministic smooth function and the law of $\bar{X}^{(1)}$ is $\delta_{\bar{X}^{(1)}_t}$, and \eqref{eq1} reduces the following non-degenerate McKean-Vlasov SDE associated with $\bar{X}^{(2)}$:
	\begin{align}\label{reduce coreq}
		\d X_t^{(2)}=l(t, \bar{X}^{(2)}_{t}, \mathscr{L}_{\bar{X}^{(2)}_{t}})\d t+\d W(t), Y(0)=y,
	\end{align}
	where $l(t, \bar{X}^{(2)}_{t}, \mathscr{L}_{\bar{X}^{(2)}_{t}}):=q\left(\bar{X}_{t}, \mathscr{L}_{\bar{X}^{(2)}_{t}}\right)$, by Theorem \ref{thm 2}, we can easily get the Onsager-Machlup action functional for \eqref{reduce coreq} is given by
	$$
	L(\phi^2, \dot{\phi}^2, \delta_{\phi^2_t})=-\frac{1}{2}\int_{0}^{1}\Big|\dot{\phi}^2(t)-l(t, \phi^{2}_{t}, \mathscr{L}_{\phi^{2}_{t}})\Big|^2\d t-\frac{1}{2} \int_{0}^{1}\operatorname{div}l(t, \phi^{2}_{t}, \mathscr{L}_{\phi^{2}_{t}})\d t,
	$$
	which coincides with our previous results of Onsager-Machlup action functional  for non-degenerate SDEs \cite{SL}.
\end{cor}
In order to apply to more realistic McKean-Vlasov models, we reformulate a equivalent regularity criteria, namely, we assume reference path $\phi \in {C}_b^2([0,T],\mathbb{R}^{d+m})$ and then let $p \in C^1(\mathbb{R}^{d+m})$ and $q \in C^{2,(1,1)}\left(\mathbb{R}^{d+m} \times \mathscr{P}_{2}\left(\mathbb{R}^{m}\right)\right)$.

\begin{thm}
	For a class of degenerate stochastic differential equation in the form of \eqref{eq1}, assume reference path $\phi \in {C}_b^2([0,T],\mathbb{R}^{d}\times \mathbb{R}^{m})$, $p \in C^1(\mathbb{R}^{d+m})$ and $q \in C^{2,(1,1)}\left(\mathbb{R}^{d+m} \times \mathscr{P}_{2}\left(\mathbb{R}^{m}\right)\right)$ and  assumptions $\mathbf{(H_0)}, \mathbf{(H_1)}, \mathbf{(H_2)},\mathbf{(H_3)}, \mathbf{(H_4)}$ hold. Then the Onsager-Machlup action functional of $\bar{X}_t$ for any norm dominating $L^2[0,T]$  exists and is given by
	$$L(\phi,\dot{\phi},\delta_{\phi^2})=-\frac{1}{2}\int_{0}^{T}\Big|\dot{\phi}^2(t)-q(\phi(t),\delta_{\phi^2(t)})\Big|^2\d t-\frac{1}{2}\int_{0}^{T}\operatorname{div}_{x^2}q(\phi(t),\delta_{\phi^2(t)})\d t,$$
	where $\operatorname{div}_{x^2}$  denote the divergence on the $\phi^2(t)\in\mathbb{R}^m$, or the following equivalent global expression:
	$$
	L(\phi, \dot{\phi},\delta_{\phi^2})=-\frac{1}{2}\int_{0}^{T}\Big|\Xi^{+}[\dot{\phi}(t)-b(\phi(t), \delta_{\phi^2(t)})]\Big|^2\d t-\frac{1}{2}\int_{0}^{T} \operatorname{div_x}[\Xi^{+}\cdot b](\phi(t), \delta_{\phi^2(t)})\d t,
	$$
	where $\operatorname{div_x}$ denote the divergence on the $\mathbb{R}^{d+m}$ and $\phi$ is the function such that
	$$\phi(t)-x^{\mu}=(\phi^1(t)-x^{(1)},\phi^2(t)-x^{(2)})\quad \text{and} \quad \d\phi^1(t)=p(\phi(t))\d t.$$
\end{thm}
Consider the following McKean-Vlasov stochastic Hamilton system:
\begin{equation}\label{hamilton Mckean}
	\left\{\begin{array}{l}
		\mathrm{d} \bar{X}_{t}^{(1)}=\bar{X}_{t}^{(2)} \mathrm{d} t,\\
		\mathrm{d} \bar{X}_{t}^{(2)}=q\left(\bar{X}_{t}, \mathscr{L}_{\bar{X}^{(2)}_{t}}\right) \mathrm{d} t+ \mathrm{d} W(t),
	\end{array}\right.
\end{equation}
\begin{cor}\label{ham}
	For a class of degenerate stochastic differential equation in the form of \eqref{hamilton Mckean}, let reference path $\phi \in {C}_b^2([0,T],\mathbb{R}^{2d})$, $p=\bar{X}^{(2)}$ and $q \in C^{2,(1,1)}\left(\mathbb{R}^{d+m} \times \mathscr{P}_{2}\left(\mathbb{R}^{m}\right)\right)$, then under assumptions $\mathbf{(H_0)}, \mathbf{(H_1)}, \mathbf{(H_2)},\mathbf{(H_3)}, \mathbf{(H_4)}$, we can write Onsager-Machlup action functional as follows:
	\begin{align}
		L(\phi, \dot{\phi},\delta_{\phi^2_t})\nonumber
		&:=E(\phi^1, \dot{\phi}^1, \ddot{\phi}^2,\delta_{\dot{\phi^1}})=-\frac{1}{2}\int_{0}^{T} \left|\ddot{\phi}^1(t)-f(\phi^1(t),\dot{\phi}^1(t), \mathscr{L}_{\dot{\phi^1}(t)})\right|^{2}\d t\nonumber\\&-\frac{1}{2}\int_{0}^{T} \operatorname{div}_{\dot{\phi}^1}{f}(\phi^1(t),\dot{\phi}^1(t), \mathscr{L}_{\dot{\phi}^1(t)})\d t,\nonumber
	\end{align}
	where $\phi(t)-x^{\mu}=(\phi^1(t)-x^{(1)},\phi^2(t)-x^{(2)})$ and $\phi^{(2)}(t)=\dot{\phi}^{(1)}(t)$.
\end{cor}
\section{Application to some models}
In this section we will determine the most probable path for some degenerate system by using OM action function. Variational principle is the bridge between OM action functional and the most probable path.
\begin{exa}
	Consider the following scalar McKean-Vlasov SDE:
	\begin{align}\label{ex2}
		\d \dot{X}_t=E(\dot{X}_t)[X_t^2-1]\d t+\d W_t, X_0=1,X_T=-1,
	\end{align}
	which can be rewriten as a system of first-order McKean-Vlasov SDEs:
	\begin{align}\label{ex22}
		\begin{aligned}
			\d\left(\begin{array}{c}
				X_t \\
				Y_t
			\end{array}\right)= \left(\begin{array}{c}
				Y_t\\
				E(Y_t)[X_t^2-1]
			\end{array}\right) \mathrm{d} t+\left(\begin{array}{c}
				0 \\
				1
			\end{array}\right) \mathrm{d} W(t).
		\end{aligned}
	\end{align}
	Noting that $1,-1$ are metastable states of system \eqref{ex2} or \eqref{ex22}. And then by Corollary \ref{ham} we can obtain the Onsager-Machlup action functional for \eqref{ex22}:
	$$L(\phi,\dot{\phi},\delta_{\phi})=-\frac{1}{2}\int_{0}^{T}\Big|\dot{\phi}^2(t)-\int_{\mathbb{R}}((\phi^1(t))^2-1)y\mu_t^{\phi^2}(\d y)\Big|^2\d t,$$
	where $\dot{\phi}^1(t)=\phi^2(t)$. Similarly, we can rewrite Onsager-Machlup action functional by  $(\phi^1,\dot{\phi^1},\ddot{\phi^1}, \delta_{\phi^1})$, i.e.
	$$L(\phi,\dot{\phi},\delta_{\phi^2})=L(\phi^1,\dot{\phi}^1,\ddot{\phi}^1,\delta_{\dot{\phi^1}})=-\frac{1}{2}\int_{0}^{T}\Big|\ddot{\phi}^1(t)-\int_{\mathbb{R}}((\phi^1(t))^2-1)y\mu_t^{\dot{\phi^1}}(\d y)\Big|^2\d t,$$
	we can determine the most probable path $\phi^1_m$ for $X_t$ by minimizing the corresponding OM action functional $L(\phi_t,\dot{\phi}_t, \delta_{\dot{\phi^1}})$ with the application of variational principle. Then recall technique in \cite{SL} Section $4$, we present Euler-Lagrange equations for $L(\phi^1,\dot{\phi^1},\ddot{\phi^1}, \delta_{\phi^1})$,
	\begin{align}
		\frac{\d^4 \phi^1}{\d t^4}+\frac{\d^2 \phi^1}{\d t^2}\int_{\mathbb{R}}2\phi^1 y\mu_t^{\dot{\phi^1}}(\d y)-2\phi^1\frac{\d \phi^1}{\d t}-\int_{\mathbb{R}}(2(\frac{\d \phi^1}{\d t})^2+\frac{\d^2 \phi^1}{\d t^2}) y\mu_t^{\dot{\phi^1}}(\d y)\nonumber\\-\int_{\mathbb{R}}2\phi^1 y\mu_t^{\phi^1}(\d y)\int_{\mathbb{R}}((\phi^1)^2-1) y\mu_t^{\dot{\phi^1}}(\d y)=0\nonumber,
	\end{align}
	with boundary conditions $(\phi^1(0),\dot{\phi}^1(0))=(1,-1)$ and $(\phi^1(T),\dot{\phi}^1(T))=(-1,1)$.
\end{exa}
\section{Proof of main results}
\subsection{Proof of Theorem \ref{thm 2}}
\begin{proof}
	We first construct auxiliary process:
	\begin{equation}\label{X^1}
		\left\{\begin{array}{l}
			\mathrm{d} \tilde{X}_t^{(1)}=p(\tilde{X}_t)  \mathrm{d} t, \\
			\mathrm{d}\tilde{X}_t^{(2)}=\dot{\phi}^2(t)\mathrm{d} t+\mathrm{d} W(t),
		\end{array}\right.
	\end{equation}
	and $\tilde{X}_t:=(\tilde{X}_t^{(1)},\tilde{X}_t^{(2)})$ with $\tilde{X}_0=x^{\mu}.$ So if we let
	$$
	\begin{aligned}
		\begin{aligned}
			&W^{\mu}(t)=W(t)-\int_{0}^{t} \zeta(s) \mathrm{d} s, \\
			&\zeta(s)=q(\tilde{X}_s, \mathscr{L}_{\tilde{X}^{(2)}_s})-\dot{\phi}^2(s), s, t \in[0, 1],
		\end{aligned}
	\end{aligned}
	$$
	then we have
	\begin{equation}
		\left\{\begin{array}{l}
			\mathrm{d} \tilde{X}_t^{(1)}=p(\tilde{X}_t) \mathrm{d} t, \quad \tilde{X}^{(1)}_0=x^{(1)}\in \mathbb{R}^d,\\
			\mathrm{d}\tilde{X}_t^{(2)}=q(\tilde{X}_t, \mathscr{L}_{\tilde{X}_t})\mathrm{d} t+\mathrm{d} W^\mu(t),\quad  X^{(2)}_0=x^{(2)}\in \mathbb{R}^m.
		\end{array}\right.\nonumber
	\end{equation}
	The Novikov condition is clearly satisfied since properties of $q$ and $\phi^2$. So by Girsanov's theorem we have  $\left\{W^\mu(t)\right\}_{t \in[0, 1]}$ is a $m$-dimensional Brownian motion on $\mathbb{R}^m$ under probability $\frac{\mathrm{d} \mathbb{Q}^\mu}{\mathrm{d} \mathbb{P}}=R^\mu$, where
	$$
	R^\mu:=\exp \left[\int_{0}^{1}\langle\zeta(t), \mathrm{d} W(t)\rangle-\frac{1}{2} \int_{0}^{1}|\zeta(t)|^{2} \mathrm{~d} t\right].
	$$
	Before we begin to estimate small probability, we first recall the structure of reference path in the sense of distribution dependent.
	$$\phi(t)-x^{\mu}=(\phi^1(t)-x^{(1)},\phi^2(t)-x^{(2)})\in \text{Cameron-Martion space } \mathcal{H}\quad \text{and} \quad \d\phi^1(t)=p(\phi(t))\d t,$$
	where $\phi^1(t)\in\mathbb{R}^{d}, \phi^2(t)\in\mathbb{R}^{m}$ and $\phi^1(0)=x^{(1)},\phi^2(t)=x^{(2)}$.  \\
	We now apply Girsanov transform to  $\mathbb{P}( \|\bar{X}-\phi\Vert<\varepsilon)$ with $\mathrm{d}W^\mu(t)=\mathrm{d}W(t)-\zeta(t)\mathrm{d}t$,  we obtain for any $\varepsilon>0.$
	\begin{align}\label{all-j1}
		&\mathbb{P}\Big(\Vert \bar{X}-\phi\Vert<\varepsilon\Big)=\mathbb{Q}^\mu\Big(\Vert \tilde{X}-\phi\Vert<\varepsilon\Big)=E\Big(R^\mu\mathbb{I}_{\Vert \tilde{X}_t-\phi\Vert<\varepsilon}\Big)\nonumber\\=&E\Big(\exp \left[\int_{0}^{1}\langle\zeta(t), \mathrm{d} W(t)\rangle-\frac{1}{2} \int_{0}^{1}|\zeta(t)|^{2} \mathrm{d} t\right]\mathbb{I}_{\Vert \tilde{X}_t-\phi\Vert<\varepsilon}\Big)\nonumber\\=&E\Big(\exp\Big[\int_{0}^{1}\langle q(\tilde{X}_t, \mathscr{L}_{\tilde{X}^{(2)}_t})-\dot{\phi}^2(t), \mathrm{d} W(t)\rangle\nonumber\\&\nonumber-\frac{1}{2}\int_{0}^{1}\Big|q(\tilde{X}_t, \mathscr{L}_{\tilde{X}^{(2)}_t})-\dot{\phi}^2(t)\Big|^2\d t\mathbb{I}_{\|\tilde{X}_t-\phi\|<\varepsilon}\Big).\nonumber\\:=&\exp(J_0) E\Big[\exp(\sum_{i=1}^{4}J_i)|\|\tilde{X}-\phi\|<\varepsilon\Big]\mathbb{P}(\|\tilde{X}-\phi\|<\varepsilon),
	\end{align}
	with
	\begin{align}
		J_0&:=\int_{0}^{1}-\frac{1}{2}|q(\phi(t),\mathscr{L}_{\phi^2(t)})-\dot\phi^2(t)|^2\mathrm{d}t,\nonumber\\
		J_1&:=\int_{0}^{1}\langle q(\tilde{X}_t, \mathscr{L}_{\tilde{X}^{(2)}_t} )-q(\phi(t), \mathscr{L}_{\phi^{(2)}(t)}),\dot\phi^2(t)\rangle \mathrm{d}t,\nonumber\\
		J_2&:=\frac{1}{2}\int_{0}^{1}|q(\phi(t), \mathscr{L}_{\phi^2(t)} )|^2\mathrm{d}t-\frac{1}{2}\int_{0}^{1}\Big|q(\tilde{X}_t, \mathscr{L}_{\tilde{X}^{(2)}_t})\Big|^2\mathrm{d}t,\nonumber\\
		J_3&:=\int_{0}^{1}\langle-\dot{\phi}^2(t),\mathrm{d}W(t)\rangle,\nonumber\\J_4&=\int_{0}^{1}\langle q(\tilde{X}_t, \mathscr{L}_{\tilde{X}^{(2)}_t}),\mathrm{d}W(t)\rangle.\nonumber
	\end{align}
	By assumption $\mathbf{(H_4)}$, Remark \ref{equiva of small ball} and \eqref{X^1}, we have
	\begin{align}\label{all-j2}
		\mathbb{P}\Big(\Vert \bar{X}-\phi\Vert<\varepsilon\Big)=\exp(J_0) E\Big[\exp(\sum_{i=1}^{4}J_i)|\|W(t)\|<\varepsilon\Big]\mathbb{P}(\|W\|<\varepsilon).
	\end{align}
	Firstly, we deal with the terms $J_1$ and $J_2$. Applying  H$\mathrm{\ddot{o}}$lder inequality and $\mathbf{(H_0)} (iii)$ we easily get
	\begin{align}\label{J_1}
		\limsup _{\varepsilon \rightarrow 0}\ E(\exp(cJ_{1})|\|W\|<\varepsilon)\leq1,
	\end{align}
	\begin{align}\label{J_2}
		\limsup _{\varepsilon \rightarrow 0}\ E(\exp(cJ_{2})|\|W\|<\varepsilon)\leq1,
	\end{align}
	for every real number $c$.\\
	Moreover, by Lemma \ref{no random lemma} we have
	\begin{align}\label{J_3}
		\limsup _{\varepsilon \rightarrow 0}\ E(\exp(cJ_{3})|\|W\|<\varepsilon)\leq1,
	\end{align}
	for every real number $c$.\\
	As for the term $J_4$, applying Lemma \ref{aa} to $q_i(\tilde{X}_t,\mathscr{L}_{\tilde{X}^{(2)}_t})W_i(t)$, we obtain
	\begin{align}
		\int_{0}^{1}q_i(\tilde{X}_t,\mathscr{L}_{\tilde{X}^{(2)}_t})\d W_i(t)&= q_i(\tilde{X}_1,\mathscr{L}_{\tilde{X}^{(2)}_1})W_i(1)-\int_{0}^{1}W_i(t)\mathbf{L}_{1}^{i}q_i(\tilde{X}^{(1)}_t,\tilde{X}^{(2)}_t,\mathscr{L}_{\tilde{X}^{(2)}_t})\d t\nonumber\\-&\sum_{j=1}^{m}\int_{0}^{1}W_i(t)\partial_{x^2_j}q_i(\tilde{X}_t,\mathscr{L}_{\tilde{X}^{(2)}_t})\d W_j(t)-\sum_{k=1}^{d}\int_{0}^{1}W_i(t)\partial_{x^1_k}q_i(\tilde{X}_t,\mathscr{L}_{\tilde{X}^{(2)}_t})\d\tilde{X}_t^{(1)}\nonumber\\-&\int_{0}^{1}\partial_{x^2_i}q_i(\tilde{X}_t,\mathscr{L}_{\tilde{X}^{(2)}_t})\d t,\nonumber
	\end{align}
	where $x=(x^1,x^2)=(x^1_1,\cdots, x^1_d,x^2_1,\cdots, x^2_m)$ and $\mathbf{L}^i_{1} q_i(x^1, x^2, \mu^{(2)})$ is defined by
	$$
	\begin{aligned}
		\mathbf{L}^i_{1} q_i(x^1, x^2, \mu^{(2)})=& \frac{1}{2} \operatorname{tr}\left( \partial_{x^2}^{2} q_i\right)(x, \mu^{(2)})+\langle \dot{\phi}^2, \partial_{x^2} q_i\rangle(t, x, \mu^{(2)}) \\
		&+\int_{\mathbb{R}^{m}}\Big[\frac{1}{2} \operatorname{tr}\left\{ \partial_{y} \partial_{\mu} q_i(x, \mu^{(2)})(y)\right\}+\langle \dot{\phi}_t, \partial_{\mu^{(2)}}q_i(x,\mu^{(2)})(y)\rangle\Big]\mu^{(2)}(\d y).
	\end{aligned}
	$$
	So we can rewrite $J_4$ as
	\begin{align}\label{part of J_4}
		J_4=& \sum_{i=1}^{m}q_i(\tilde{X}_1,\mathscr{L}_{\tilde{X}^{(2)}_1})W_i(1)-\sum_{i=1}^{m}\int_{0}^{1}W_i(t)\mathbf{L}_{1}^{i}q_i(\tilde{X}^{(1)}_t,\tilde{X}^{(2)}_t,\mathscr{L}_{\tilde{X}^{(2)}_t})\d t\nonumber\\-&\sum_{i=1}^{m}\sum_{j=1}^{m}\int_{0}^{1}W_i(t)\partial_{x^2_j}q_i(\tilde{X}_t,\mathscr{L}_{\tilde{X}^{(2)}_t})\d W_j(t)-\sum_{i=1}^{m}\sum_{k=1}^{d}\int_{0}^{1}W_i(t)\partial_{x^1_k}q_i(\tilde{X}_t,\mathscr{L}_{\tilde{X}^{(2)}_t})\d\tilde{X}_t^{(1)}\nonumber\\-&\sum_{i=1}^{m}\int_{0}^{1}\partial_{x^2_i}q_i(\tilde{X}_t,\mathscr{L}_{\tilde{X}^{(2)}_t})\d t,\nonumber\nonumber\\:=&J_5+J_6+J_7+J_8,
	\end{align}
	where
	$$
	\begin{aligned}
		J_5=&\sum_{i=1}^{m}q_i(\tilde{X}_1,\mathscr{L}_{\tilde{X}^{(2)}_1})W_i(1),\nonumber\\J_6=&-\sum_{i=1}^{m}\int_{0}^{1}W_i(t)\mathbf{L}_{1}^{i}q_i(\tilde{X}^{(1)}_t,\tilde{X}^{(2)}_t,\mathscr{L}_{\tilde{X}^{(2)}_t})\d t-\sum_{i=1}^{m}\sum_{k=1}^{d}\int_{0}^{1}W_i(t)\partial_{x^1_k}q_i(\tilde{X}_t,\mathscr{L}_{\tilde{X}^{(2)}_t})\d\tilde{X}_t^{(1)},\nonumber\\J_7=&-\sum_{i=1}^{m}\sum_{j=1}^{m}\int_{0}^{1}W_i(t)\partial_{x^2_j}q_i(\phi(t),\mathscr{L}_{\phi^2(t)})\d W_j(t)-\sum_{i=1}^{m}\int_{0}^{1}\partial_{x^2_i}q_i(\tilde{X}_t,\mathscr{L}_{\tilde{X}^{(2)}_t})\d t,\\ J_8=&-\sum_{i=1}^{m}\sum_{j=1}^{m}\int_{0}^{1}W_i(t)\Big[\partial_{x^2_j}q_i(\tilde{X}_t,\mathscr{L}_{\tilde{X}^{(2)}_t})-\partial_{x^2_j}q_i(\phi(t),\mathscr{L}_{\phi^2(t)})\Big]\d W_j(t).
	\end{aligned}
	$$
	For the term $J_5$, as $q_i (i=1,\cdots, m)$ are bounded, by Lemma \ref{Separation lemma} and \ref{no random lemma}, we have
	\begin{align}\label{J_5}
		\limsup _{\varepsilon \rightarrow 0}\ E(\exp(cJ_{5})|\|W\|<\varepsilon)\leq1,
	\end{align}
	for every real number $c$.\\
	Then we deal with the term $J_6$. Recall the expression of $d\tilde{X}^{(1)}_t=p(\tilde{X}_t) \mathrm{d} t$, so we get
	\begin{align}
		J_6=&-\sum_{i=1}^{m}\int_{0}^{1}W_i(t)\mathbf{L}_{1}^{i}q_i(\tilde{X}^{(1)}_t,\tilde{X}^{(2)}_t,\mathscr{L}_{\tilde{X}^{(2)}_t})\d t-\sum_{i=1}^{m}\sum_{k=1}^{d}\int_{0}^{1}W_i(t)\partial_{x^1_k}q_i(\tilde{X}_t,\mathscr{L}_{\tilde{X}^{(2)}_t})\d\tilde{X}_t^{(1)}\nonumber\\=&-\sum_{i=1}^{m}\int_{0}^{1}W_i(t)\mathbf{L}_{1}^{i}q_i(\tilde{X}^{(1)}_t,\tilde{X}^{(2)}_t,\mathscr{L}_{\tilde{X}^{(2)}_t})\d t-\sum_{i=1}^{m}\sum_{k=1}^{d}\int_{0}^{1}W_i(t)\partial_{x^1_k}q_i(\tilde{X}_t,\mathscr{L}_{\tilde{X}_t})p(\tilde{X}_t, \mathscr{L}_{\tilde{X}^{(2)}_t}) \mathrm{d} t\nonumber
	\end{align}
	Since $q \in C_{b}^{2,(1,1)}\left( \mathbb{R}^{d+m} \times \mathscr{P}_{2}\left(\mathbb{R}^{m}\right)\right)$ and $p \in C_{b}^{1}\left( \mathbb{R}^{d+m}\right)$,  we get for $1\leq i \leq m, 1\leq k \leq d$
	\begin{align}
		|\int_{0}^{1}W_i(t)\mathbf{L}_{1}^{i}q_i(\tilde{X}^{(1)}_t,\tilde{X}^{(2)}_t,\mathscr{L}_{\tilde{X}^{(2)}_t})dt|\leq K_5|\int_{0}^{1}W_i(t)\d t|\leq K_5\varepsilon,\nonumber\\|\int_{0}^{1}W_i(t)\partial_{x^1_k}q_i(\tilde{X}_t,\mathscr{L}_{\tilde{X}^{(2)}_t})p(\tilde{X}_t) \mathrm{d}t|\leq K_6|\int_{0}^{1}W_i(t)\d t|\leq K_6\varepsilon.\nonumber
	\end{align}
	By using Lemma \ref{Separation lemma} and \ref{no random lemma} again, we obtain
	\begin{align}\label{J_6}
		\limsup _{\varepsilon \rightarrow 0}\ E(\exp(cJ_{6})|\|W\|<\varepsilon)\leq1,
	\end{align}
	for every real number $c$.\\
	For the term $J_7$, we first divide it into four parts as usual.
	\begin{align}
		J_7=&-\sum_{i\neq j}\int_{0}^{1}W_i(t)\partial_{x^2_j}q_i(\phi(t),\mathscr{L}_{\phi^2(t)})\d W_j(t)-\sum_{i=1}^{m}\int_{0}^{1}W_i(t)\partial_{x^2_i}q_i(\phi(t),\mathscr{L}_{\phi^2(t)})\d W_i(t)\nonumber\\&-\sum_{i=1}^{m}\int_{0}^{1}\Big[\partial_{x^2_i}q_i(\tilde{X}_t,\mathscr{L}_{\tilde{X}^{(2)}_t})-\partial_{x^2_i}q_i(\phi(t),\mathscr{L}_{\phi^2(t)})\Big]\d t-\sum_{i=1}^{m}\int_{0}^{1}\partial_{x^2_i}q_i(\phi(t),\mathscr{L}_{\phi^2(t)})\d t.\nonumber
	\end{align}
	By Lemma \ref{Cross lemma}, when $i\neq j$, we have
	\begin{align}
		\limsup _{\varepsilon \rightarrow 0}\ E(\exp(c-\sum_{i\neq j}\int_{0}^{1}W_i(t)\partial_{x^2_j}q_i(\phi(t),\mathscr{L}_{\phi^2(t)})\d W_j(t))|\|W\|<\varepsilon)\leq1,\nonumber
	\end{align}
	for every real number $c$.\\
	Applying Lemma \ref{aa} to $\partial_{x^2_i}q_i(\phi(t),\mathscr{L}_{\phi^2(t)})W^2_i(t)$, we obtain
	\begin{align}
		\int_{0}^{1}\partial_{x^2_i}q_i(\phi(t),\mathscr{L}_{\phi^2(t)})W_i(t)\d W_i(t)&= \frac{1}{2}\partial_{x^2_i}q_i(\phi(1),\mathscr{L}_{\phi^2(1)})(W_i(1))^2\nonumber\\&-\frac{1}{2}\int_{0}^{1}(W_i(t))^2\mathbf{L}_{2}^{i}q_i(\phi(t),\mathscr{L}_{\phi^2(t)})\d t-\frac{1}{2}\int_{0}^{1}\partial_{x^2_i}q_i(\phi(t),\mathscr{L}_{\phi^2(t)})\d t,
	\end{align}
	where
	\begin{align}
		\mathbf{L}_2^{i}\partial_{x_i}q_i(x,\mu^{(2)})= \langle \dot{\phi}, \partial_{x} \partial_{x^2_i}q_i\rangle(x, \mu^{(2)}) +\int_{\mathbb{R}^{m}}\langle\dot{\phi}_t, \partial_{\mu^{(2)}}[\partial_{x_i}f_i(t,x,\mu^{(2)})](y)\rangle\mu^{(2)}(\d y).\nonumber
	\end{align}
	Since the derivatives of $q_i$ are bounded, so by
	$\mathbf{(H_2)}$, we have
	\begin{align}
		\limsup _{\varepsilon \rightarrow 0}\ E(\exp(\frac{c}{2}\partial_{x^2_i}q_i(\phi(1),\mathscr{L}_{\phi^{(2)}(1)})(W_i(1))^2)|\|W\|<\varepsilon)\leq1,\nonumber
	\end{align}
	for  every real number $c$.\\
	By the properties of $q \in C_{b}^{2,(1,1)}\left( \mathbb{R}^{d+m} \times \mathscr{P}_{2}\left(\mathbb{R}^{m}\right)\right)$, we have
	\begin{align}
		|-\frac{1}{2}\int_{0}^{1}(W_i(t))^2\mathbf{L}_{2}^{i}q_i(\phi(t),\mathscr{L}_{\phi^{(2)}(t)})\d t|\leq K_7\varepsilon^2.\nonumber
	\end{align}
	As a consequence,
	\begin{align}
		\limsup _{\varepsilon \rightarrow 0}\ E(-\frac{c}{2}\int_{0}^{1}(W_i(t))^2\mathbf{L}_{2}^{i}q_i(\phi(t),\mathscr{L}_{\phi^{(2)}(t)})\d t)|\|W\|<\varepsilon)\leq1,\nonumber
	\end{align}
	for  every real number $c$.\\
	By $\mathbf{(H_0)} (iii)$ and Lemma \ref{Separation lemma}, we have that
	\begin{align}
		\limsup _{\varepsilon \rightarrow 0}\ E(\exp(c|\sum_{i=1}^{m}\int_{0}^{1}\Big[\partial_{x^2_i}q_i(\tilde{X}_t,\mathscr{L}_{\tilde{X}^{(2)}_t})-\partial_{x^2_i}q_i(\phi(t),\mathscr{L}_{\phi^2(t)})\Big]\d t|\|W\|<\varepsilon)\leq1,\nonumber
	\end{align}
	for  every real number $c$.\\
	Therefore by Lemma \ref{Separation lemma}, we obtain
	\begin{align}\label{J_7}
		\limsup _{\varepsilon \rightarrow 0}\ E(\exp(c|J_7+\frac{1}{2}\int_{0}^{1}\operatorname{div}_{x^2}{q}(\phi(t), \mathscr{L}_{\phi^2(t)} )\d t|\|W\|<\varepsilon)\leq1,
	\end{align}
	for  every real number $c$. And $-\frac{1}{2}\int_{0}^{1}\operatorname{div}_{x^2}{q}(\phi(t),\mathscr{L}_{\phi^2(t)})\d t$ is the divergence part of OM action functional for degenerate McKean-Vlasov SDEs.\\
	Finally, we deal with term $J_{8}$ by estimation of small ball probability. We first define $T_{8}^{ij}:=\int_{0}^{1}cW_i(t)\Big[\partial_{x^2_j}q_i(\tilde{X}_t,\mathscr{L}_{\tilde{X}^{(2)}_t})-\partial_{x^2_j}q_i(\phi(t),\mathscr{L}_{\phi^2(t)})\Big]\d W_j(t)$, and for any $c\in \mathbb{R}$ and $\delta>0$ we get
	\begin{align}
		E&\Big(\exp(J^{ij}_{8})|\|W\|\leq \varepsilon\Big)\nonumber\\&\leq e^{\delta}+\int_{\delta}^{\infty}e^{\xi}\mathbb{P}(\Big|J^{ij}_8\Big|>\xi|\|W\|\leq \varepsilon)\d \xi \nonumber\\&+e^{\delta}\mathbb{P}(\Big|J^{ij}_8|>\delta \Big|\|W\|\leq \varepsilon).\nonumber
	\end{align}
	Since $\int_{0}^{t}cW_i(s)\Big[\partial_{x^2_j}q_i(\tilde{X}_s,\mathscr{L}_{\tilde{X}^{(2)}_s})-\partial_{x^2_j}q_i(\phi(s),\mathscr{L}_{\phi^2(s)})\Big]\d W_j(s)$ is a martingale, and its quadratic variations can be estimated by $\mathbf{(H_0)}$ (iii), $(\mathbf{H_3})$, Young's inequality as follows,
	\begin{align}
		\langle &\int_{0}^{t}cW_i(s)\Big[\partial_{x^2_j}q_i(\tilde{X}_s,\mathscr{L}_{\tilde{X}^{(2)}_s})-\partial_{x^2_j}q_i(\phi(s),\mathscr{L}_{\phi^2(s)})\Big]\d W_j(s) \rangle_t\nonumber\\\nonumber=&\int_{0}^{t}c^2\Big|W_i(s)\Big(\partial_{x^2_j}q_i(\tilde{X}_s,\mathscr{L}_{\tilde{X}^{(2)}_s})-\partial_{x^2_j}q_i(\phi(s),\mathscr{L}_{\phi^2(s)})\Big)\Big|^2\d s\nonumber\\&\leq K_{8}\Big[\int_{0}^{1}|\tilde{X}^{(1)}_t|^4\d t+\int_{0}^{1}|W(t)|^4\d t\Big]\nonumber\\&\leq K_{9}\varepsilon^{p\wedge 4},\nonumber
	\end{align}
	for $1\leq i,j \leq m$. So by the standard exponential inequality for martingales, we have
	\begin{align}
		\mathbb{P}(|T^{ij}_{8}|>\delta, \|W\|\leq \varepsilon)\leq \exp\Big(-\frac{\delta^2}{2K_{9}\varepsilon^{p\wedge 4}}\Big).\nonumber
	\end{align}
	Combining assumption $\mathbf{(H_3)}$, we obtain
	\begin{align}
		\mathbb{P}(\Big|cJ_8^{ij} \Big|>\xi|\|W\|<\varepsilon)&=\frac{\mathbb{P}(\Big|cJ_8^{ij}\Big|>\xi, \|W\|\leq \varepsilon)}{\mathbb{P}(\|W\|\leq\varepsilon)}\nonumber\\&\leq \exp\Big\{-\frac{\xi^2}{2K_{9}\varepsilon^{p\wedge 4}}+\frac{C_3}{\varepsilon^q}\Big\},\nonumber
	\end{align}
	for every real number $c$.\\
	Using the previous estimates we have for every $\delta>0$ and every $0<\varepsilon<1$
	\begin{align}
		E&\Big(\exp(cJ^{ij}_{8})|\|W\|\leq \varepsilon\Big)\nonumber\\&\leq e^{\delta}+\int_{\delta}^{\infty}\exp\Big\{\xi-\frac{\xi^2}{2K_{9}\varepsilon^{p\wedge 4}}+\frac{C_3}{\varepsilon^q}\Big\}\d \xi \nonumber\\&+\exp\Big\{\delta-\frac{\delta^2}{2K_{9}\varepsilon^{p\wedge 4}}+\frac{C_3}{\varepsilon^q}\Big\}.\nonumber
	\end{align}
	Letting $\varepsilon$ and then $\delta$ tend to zero, we obtain
	\begin{align}
		\limsup _{\varepsilon \rightarrow 0}\ E(\exp(J^{ij}_{8})|\|W\|<\varepsilon)\leq1,\nonumber
	\end{align}
	for every real number $c$ and $1\leq i,j \leq m$. As a consequence, by Lemma \ref{Separation lemma} we have
	\begin{align}\label{J_8}
		\limsup _{\varepsilon \rightarrow 0}\ E(\exp(cJ_{8})|\|W\|<\varepsilon)\leq1,
	\end{align}
	for every real number $c$.\\
	Recall previous estimates \eqref{all-j1},\eqref{all-j2},\eqref{J_1}, \eqref{J_2}, \eqref{J_3},\eqref{part of J_4}, \eqref{J_5}, \eqref{J_6}, \eqref{J_7}, \eqref{J_8} and then Lemma \ref{Separation lemma} allows us to conclude that
	$$
	\lim _{\varepsilon \rightarrow 0} \frac{\mathbb{P}(\|\bar{X}-\phi\|<\varepsilon)}{\mathbb{P}(\|W\|<\varepsilon)}=\exp \left(-\frac{1}{2} \int_{0}^{1}|\dot{\phi}^2(t)-q(\phi(t),\delta_{\phi^2(t)})|^{2} \d t-\frac{1}{2} \int_{0}^{1}\operatorname{div}_{x^2}q(\phi(t),\delta_{\phi^2(t)})\d t\right).
	$$
	The proof of Theorem \ref{thm 2} is complete.
\end{proof}


\beg{thebibliography}{99}

\bibitem{SA}
S. Aihara and A. Bagchi, On the {M}ortensen equation for maximum likelihood state
estimation, {\it IEEE Trans. Automat. Control}, {\bf 44}, $1999$, 1955-1961.

\bibitem{BA1}
B. Ayanbayev, I. Klebanov, H. C. Lie, T. J. Sullivan, $\Gamma$-convergence of Onsager-Machlup functionals. Part I: With applications to maximum a posteriori estimation in Bayesian inverse problems, {\it Inverse Problems}, {\bf 38}, $2022$, Paper No. 025005, 32.

\bibitem{XB1}
X. Bardina, C. Rovira and S. Tindel, Onsager-{M}achlup functional for stochastic evolution
equations, {\it Ann. Inst. H. Poincar\'{e} Probab. Statist.}, {\bf 39}, $2003$, 69-93.

\bibitem{AB}
A. Ben-Israel and T. N. E. Greville, Generalized inverses: theory and applications, Robert E. Krieger Publishing Co., Inc., Huntington, N.Y., $1974$.

\bibitem{GB}
G. Bet, V. Jacquier and F. R. Nardi, Effect of energy degeneracy on the transition time for a
series of metastable states, {\it J. Stat. Phys.}, {\bf 184}, $2021$, Paper No. 8, 42.

\bibitem{CG}
M. Carfagnini and  M. Gordina, On the Onsager-Machlup functional for the Brownian motion on the Heisenberg group, {\it arxiv: 1908.09182}, $2023$.

\bibitem{MC1}
M. Capitaine, Onsager-{M}achlup functional for some smooth norms on {W}iener
space, {\it Probab. Theory Related Fields}, {\bf 102}, $1995$, 189-201.

\bibitem{MC2}
M. Capitaine, On the {O}nsager-{M}achlup functional for elliptic diffusion
processes, In S\'{e}minaire de {P}robabilit\'{e}s, {XXXIV} {\it Lecture Notes in Math.}, {\bf 1729}, $2000$, 313-328.

\bibitem{YC}
Y. Chao and J. Duan,  The {O}nsager-{M}achlup function as {L}agrangian for the most probable path of a jump-diffusion process, {\it Nonlinearity}, {\bf 32}, $2019$, 3715-3741.

\bibitem{MCM}
M. Chaleyat-Maurel and D. Nualart,  Onsager-{M}achlup functionals for solutions of stochastic
boundary value problems, In S\'{e}minaire de {P}robabilit\'{e}s, {XXIX}, Springer, Berlin, Heidelberg, $1995$, 44-55.

\bibitem{MD}
M. Dashti, K. J. H. Law, A. M. Stuart and J. Voss, M{AP} estimators and their consistency in {B}ayesian
nonparametric inverse problems, {\it Inverse Problems}, {\bf 29}, $2013$, 095017, 27.

\bibitem{QD}
Q. Du, T. Li, X. Li and W. Ren, The graph limit of the minimizer of the {O}nsager-{M}achlup
functional and its computation, {\it Sci. China Math.}, {\bf 64}, $2021$, 239-280.

\bibitem{DB}
D. D\"{u}rr and A. Bach, The {O}nsager-{M}achlup function as {L}agrangian for the most probable path of a diffusion process, {\it Comm. Math. Phys.}, {\bf 60}, $1978$, 153-170.

\bibitem{AD}
A. Dembo and O. Zeitouni, Onsager-{M}achlup functionals and maximum a posteriori
estimation for a class of non-{G}aussian random fields, {\it J. Multivariate Anal.}, {\bf 36}, $1991$, 243-262.

\bibitem{WE}
W. E and E. Vanden-Eijnden, Towards a theory of transition paths, {\it J. Stat. Phys.}, {\bf 123}, $2006$, 503--523.

\bibitem{TF}
T. Fujita and S.-I. Kotani, The {O}nsager-{M}achlup function for diffusion processes, {\it J. Math. Kyoto Univ.}, {\bf 22}, $1982$, 115-130.

\bibitem{FW}
M.I. Freidlin and A.D. Wentzell, Random perturbations of dynamical systems, Springer-Verlag, New York, 1984.

\bibitem{GS}
E. Grong and S. Sommer, Most probable paths for anisotropic Brownian motions on manifolds, {\it arxiv: 2110.15634}, $2022$.

\bibitem{KH}
K. Hara and Y. Takahashi, Lagrangian for pinned diffusion process,  Springer, {\it It\^{o}'s stochastic calculus and probability theory},  $1996$, 117-128.

\bibitem{KH1}
K. Hara and Y. Takahashi, Stochastic analysis in a tubular neighborhood or Onsager-Machlup functions revisited, {\it arXiv:1610.06670}, $2016$.

\bibitem{JH}

J. Hu, X. Chen and J. Duan, An Onsager-Machlup approach to the most probable transition
pathway for a genetic regulatory network, {\it Chaos}, {\bf 32}, $2022$, Paper No. 041103.

\bibitem{YH1}
Y. Hu, Multi-dimensional geometric {B}rownian motions,
{O}nsager-{M}achlup functions, and applications to
mathematical finance, {\it Acta Math. Sci. Ser. B}, {\bf 20}, $2000$, 341-358.

\bibitem{YH2}
Y. Huang, Y. Chao and W. Wei, Estimating the most probable transition time for stochastic dynamical systems, {\it Nonlinearity}, {\bf 34}, $2021$, 4543--4569.

\bibitem{HH}
C. H. Hung and T. L. Markham, The {M}oore-{P}enrose inverse of a partitioned matrix
$M=\left(\begin{array}{ll}
	A & B \\
	C & D
\end{array}\right)$, {\it Linear Algebra Appl.}, {\bf 11}, $1975$, 73--86.

\bibitem{NI}
N. Ikeda and S. Watanabe, Stochastic differential equations and diffusion processes, Elsevier, $2014$.

\bibitem{JK}
J. Kurchan, Fluctuation theorem for stochastic dynamics, {\it J. Phys. A}, {\bf 31}, $1998$, 3719--3729.

\bibitem{TL}
T. Li and X. Li, Gamma-limit of the {O}nsager-{M}achlup functional on the space
of curves,  {\it SIAM J. Math. Anal.}, {\bf 53}, $2021$, 1-31.

\bibitem{QL}
Q. Liu, D. Jiang,  T. Hayat and A. Alsaedi, Dynamics of a stochastic {SIR} epidemic model with distributed
delay and degenerate diffusion, {\it J. Franklin Inst.}, {\bf 356}, $2019$, 7347--7370.

\bibitem{SL}
S. Liu, H. Gao,  H. Qiao and N. Lu, The Onsager-Machlup action functional for McKean-Vlasov SDEs, {\it Commun. Nonlinear Sci. Numer. Simul.}, 121(2023), 15, 107203.

\bibitem{WL}
W. Liu, Y. Song and T. Zhang, Large and moderate deviation principles for McKean-Vlasov SDEs with jumps, {\it arXiv:2011.08403}, $2020$.

\bibitem{SM}
S. Moret and D. Nualart,  Onsager-Machlup functional for the fractional Brownian
motion, {\it Probab. Theory Related Fields}, {\bf 124}, $2002$, 227-260.

\bibitem{OM1}
L. Onsager and S. Machlup, Fluctuations and irreversible processes, {I}, {\it Phys. Rev.}, {\bf 91}, $1953$, 1505-1512.

\bibitem{OM2}
L. Onsager and S. Machlup, Fluctuations and irreversible processes, {II}, {\it Phys. Rev.}, {\bf 91}, $1953$, 1512-1515.

\bibitem{PR}
P. Ren and F. Wang, Space-distribution {PDE}s for path independent additive
functionals of {M}c{K}ean-{V}lasov {SDE}s,  {\it Infin. Dimens. Anal. Quantum Probab. Relat. Top.}, {\bf 23}, $2020$, 2050018, 15.

\bibitem{HR}
H. Risken, The {F}okker-{P}lanck equation: Methods of solution and applications, Springer-Verlag, Berlin, $1989$.

\bibitem{RL}
R. L. Stratonovich, On the probability functional of diffusion processes, {\it Sel. Trans. Math. Stat. Prob.}, {\bf 10}, $1957$, 273-286.

\bibitem{LS}
L. Shepp and O. Zeitouni, A note on conditional exponential moments and
{O}nsager-{M}achlup functionals, {\it Ann. Probab.}, {\bf 20}, $1992$, 652-654.

\bibitem{LS1}
L. Shepp and O. Zeitouni, Exponential estimates for convex norms and some applications,  Barcelona {S}eminar on {S}tochastic {A}nalysis ({S}t. {F}eliu
de {G}u\'{\i}xols, 1991), Progr. Probab, {\bf 32}, $1993$, 203--215.

\bibitem{ZS}
Z. Selk, W. Haskell and H. Honnappa, Information projection on {B}anach spaces with applications to
state independent {KL}-weighted optimal control, {\it Appl. Math. Optim.}, {\bf 84}, $2021$, S805--S835.

\bibitem{YT}
Y. Takahashi and S. Watanabe, The probability functionals ({O}nsager-{M}achlup functions) of diffusion processes, In Stochastic integrals, Springer, $1981$, 433-463.

\bibitem{TT1}
T. Taniguchi and E. G. D. Cohen, Onsager-{M}achlup theory for nonequilibrium steady states and
fluctuation theorems, {\it J. Stat. Phys.}, {\bf 126}, $2007$, 1--41.

\bibitem{TT2}
T. Taniguchi and E. G. D. Cohen, Inertial effects in nonequilibrium work fluctuations by a path
integral approach, {\it J. Stat. Phys.}, {\bf 130}, $2008$, 1--26.

\bibitem{LT}
L. Tisza and I. Manning, Fluctuations and irreversible thermodynamics, {\it Phys. Rev.}, {\bf 105}, $1957$, 1695-1705.

\bibitem{AT}
A. Tesfay, S. Yuan, D. Tesfay and J. Brannan, Most Probable Dynamics of the Single-Species with Allee Effect under
Jump-diffusion Noise, {\it arXiv:2112.07234}, $2021$.

\bibitem{PW}
P. Wang and G. Chen, Invariant behavior of stochastic atmosphere-ocean model with degenerate noise,  {\it J. Math. Phys.}, {\bf 60}, $2019$, 062701, 17pp.

\bibitem{ZW}
Z. Wang and X. Zhang, Existence and uniqueness of degenerate {SDE}s with {H}\"{o}lder
diffusion and measurable drift, {\it J. Math. Anal. Appl.}, {\bf 484}, $2020$, 123679, 11pp.

\bibitem{OZ}
O. Zeitouni and A. Dembo, A maximum a posteriori estimator for trajectories of diffusion
processes, {\it Stochastics}, {\bf 20}, $1987$, 221-246.

\bibitem{OZ1}
O. Zeitouni and A. Dembo, An existence theorem and some properties of maximum a
posteriori estimators of trajectories of diffusions, {\it Stochastics}, {\bf 23}, $1988$, 197--218.

\bibitem{YZ}
Y. Zheng, F. Yang, J. Duan, X. Sun, L. Fu and J. Kurth, The maximum likelihood climate change for global warming under
the influence of greenhouse effect and {L}\'{e}vy noise, {\it Chaos.}, {\bf 30}, $2020$, 013132.
\end{thebibliography}

\section{Appendix}
\subsection{Verify assumption $\mathbf{(H_4)}$ by some specific norms}
Here we first verify the reasonable of hypothesis $(\mathbf{H_4})$ by  $L^p$ and H$\mathrm{\ddot{o}}$lder-norms.

$\bullet$ $L^p$-norm.
Firstly, by H$\mathrm{\ddot{o}}$lder inequality, $|a+b|^r\leq \max\{1, 2^{r-1}\}(|a|^r+|b|^r)$ (for $r>0$) and Lipschitz continuous of $g$, there exists a constant $\tau_1$ such that
\begin{align}
|\bar{X}^{(1)}_t-\phi^1(t)|^p=&\Big|\int_{0}^{t}[g(\bar{X}_t^{(1)}, \bar{X}_t^{(2)})-g(\phi^1(s), \phi^2(s))]\d s\Big|^p\leq \int_{0}^{t}|g(\bar{X}_t^{(1)}, \bar{X}_t^{(2)})-g(\phi^1(s), \phi^2(s))|^p \d s\nonumber\\\nonumber &\leq \tau_1 \int_{0}^{t}\Big[|\bar{X}_t^{(1)}-\phi^1(s)|^p+|\bar{X}_t^{(2)}-\phi^2(s)|^p\Big]\d s,
\end{align}
then by Gronwall's inequality we have for $t\in [0,1]$
\begin{align}
|\bar{X}_t^{(1)}-\phi^1(t)|^p\leq \exp\{\tau_1t\}\tau_1\int_{0}^{t}|\bar{X}_t^{(2)}-\phi^2(s)|^p\d s.
\end{align}
Therefore, by $(6.1)$ we obtain
\begin{align}
\|\bar{X}_t^{(1)}-\phi^1(t)\|^p_{L^p([0,1],\mathbb{R}^d)}\leq \tau_2\int_{0}^{1}\int_{0}^{t}|\bar{X}_s^{(2)}-\phi^2(s)|^p\d s\d t\leq \tau_2\|\bar{X}_t^{(2)}-\phi^2(t)\|^p_{L^p([0,1],\mathbb{R}^m)}.
\end{align}
So the assumption $(\mathbf{H_4})$  is valid for $L^p$-norm.

$\bullet$ H$\mathrm{\ddot{o}}$lder-norm.
In fact, we only need to verify H$\mathrm{\ddot{o}}$lder seminorm. Since supremum norm has verified in \cite{SA}.
When $p=1$, replacing $\bar{X}_t^{(1)}, \phi^1(t)$ by $\bar{X}_t^{(1)}-\phi^1(t), \bar{X}_s^{(2)}-\phi^1(s)$ $(s<t)$ respectively, then $(6.1)$ can be  rewrited as
\begin{align}
|\bar{X}_t^{(1)}-\phi^1(t)-(\bar{X}_s^{(1)}-\phi^1(s))|\leq \exp\{\tau_1(t-s)\}\tau_1\int_{s}^{t}|\bar{X}_u^{(2)}-\phi^2(u)|\d u.
\end{align}
By the definition of H$\mathrm{\ddot{o}}$lder seminorm, we obtain
\begin{align}
\|\bar{X}_t^{(1)}-\phi^1(t)\|_{C^{\alpha}([0,1],\mathbb{R}^d)}\leq \exp\{\tau_1(t-s)\}\frac{\tau_1}{|t-s|^{\alpha}}\int_{s}^{t}|\bar{X}_u^{(2)}-\phi^2(u)|\d u.
\end{align}
Let $v=\frac{u-s}{t-s}$, we have
\begin{align}
\|\bar{X}_t^{(1)}-\phi^1(t)\|_{C^{\alpha}([0,1],\mathbb{R}^d)}&\leq\exp\{\tau_1(t-s)\}\frac{\tau_1}{|t-s|^{\alpha}}\int_{0}^{1}|\bar{X}^{(2)}_{(t-s)v+s}-\phi^2((t-s)v+s)||t-s|\d v
\nonumber\\&\leq \tau_3|t-s|^{1-\alpha}\int_{0}^{1}|\bar{X}^{(2)}_{(t-s)v+s}-\phi^2((t-s)v+s)|\d v\nonumber\\&\leq \tau_4 \|\bar{X}^{(2)}_{t}-\phi^2(t)\|_{C^{\alpha}([0,1],\mathbb{R}^m)}\nonumber.
\end{align}
So the  assumption $\mathbf{H_4}$  is also valid for H$\mathrm{\ddot{o}}$lder-norm.

\end{document}